\documentclass[a4paper]{article}
\usepackage{RR}
\usepackage{hyperref}
\RRdate{D\'ecembre 2008}
\usepackage{float,amsmath,amssymb,color,epsf,graphicx}  
\RRauthor{Siddhartha Mishra \thanks{Center of Mathematics  for Applications,
University of Oslo, Oslo-0316, Norway (siddharm@cma.uio.no)} \and 
J\'er\^ome Jaffr\'e  \thanks{INRIA, BP 105, 78153 Le
Chesnay Cedex, France (Jerome.Jaffre@inria.fr)}}
\authorhead{Siddhartha Mishra and J\'er\^ome Jaffr\'e}
\RRetitle{On the upstream mobility scheme for two-phase flow\\ in porous media}
\RRtitle{Sur le sch\'ema aux mobilit\'es amont pour les \'ecoulements diphasiques\\ en milieu poreux}
\titlehead{On the upstream mobility scheme}

\newtheorem{theorem} {{\bf THEOREM}}[section]

\newtheorem{lemma} {{\bf LEMMA}}[section]

\newcommand{\be} {\begin{equation}}
\newcommand{\ee} {\end{equation}}
\newcommand{\bea} {\begin{eqnarray}}
\newcommand{\eea} {\end{eqnarray}}
\newcommand{\Bea} {\begin{eqnarray*}}
\newcommand{\Eea} {\end{eqnarray*}}
\newcommand{\dsp}{\displaystyle}
\newcommand{\pa} {\partial}

\newcommand{\De} {\Delta}
\newcommand{\la} {\lambda}

\newcommand{\no} {\nonumber}
\newcommand{\noi} {\noindent}

\newcommand{\ov} {\overline}
\newcommand{\un} {\underline}

\newcommand{\Z}{\mathbb Z}
\newcommand{\R}{\mathbb R}

\newcommand{\p}{\partial}
\newcommand{\cqfd}{\hfill \rule{1.5mm}{3mm}}

\numberwithin{equation}{section}
\numberwithin{figure}{section}

\RRabstract{
When neglecting capillarity, two-phase incompressible flow in porous
media is modelled as a scalar nonlinear hyperbolic conservation law.
A change in the rock type results in a change of the flux
function. Discretizing in one-dimensional with a finite volume method,
we investigate two numerical fluxes, an extension of the Godunov flux
and the upstream mobility flux, the latter being widely used in hydrogeology
and petroleum engineering. Then, in the case of a changing rock type,
one can give examples when the upstream mobility flux does not give
the right answer.
}

\RRresume{
En n\'egligeant la capillarit\'e, un \'ecoulement diphasique incompressible est mod\'elis\'e par une loi de conservation scalaire hyperbolique non-lin\'eaire. Un changement dans le type de roche entra\^{\i}ne un changement de la fonction flux. En discr\'etisant en dimension un avec une m\'ethode de volumes finis nous \'etudions deux flux num\'eriques, une extension du flux de Godunov et le flux des mobilt\'es amont, ce dernier \'etant largement utilis\'e en hydrologie et en ing\'eni\'erie p\'etroli\`ere. Dans le cas d'un changement de type de roche on peut alors donner des exemples o\`u le flux des mobilt\'es amont ne donne pas une solution correcte.
}

\RRmotcle{Ecoulement diphasique en milieu poreux, mobilit\'es amont, lois de conservation hyperboliques, condition d'entropie, m\'ethode de diff\'erences finies, m\'ethode de volumes finis.}
\RRkeyword{Two-phase flow in porous media, upstream mobilities, hyperbolic conservation laws, entropy condition, finite difference method, finite volume method.}

\RRprojets{Estime}
\RRtheme{\THNum}
\RCParis

\begin{document}
\RRNo{6789}
\makeRR

\section{Introduction}
Under the assumptions that capillary effects are neglected, two-phase
flow in a porous medium is modeled by a nonlinear hyperbolic
equation. In many applications, the porous medium is not
homogenous. The flow domain has to be divided into several subdomains
corresponding to different types of rock separated by lines or
surfaces along which, not only the porosity and the absolute
permeability of the rock type change but the relative permeabilities
also differ. This situation is modeled by a single conservation law
with a flux function discontinuous in the space variable. Numerical
methods designed to simulate the flow have to be devised to take into
account the discontinuities in the flux function.  

In this paper, we compare different numerical schemes of the finite
difference or finite volume type that are used to simulate two-phase
flow in porous media with changing rock types. We restrict ourselves
to the one dimensional case. In the multidimensional case, most
numerical methods still use the one dimensional flux calculation in
the direction normal to the boundaries of the discretization cells.
We also focus on the numerical flux calculation. 

Conservation laws with discontinuous coefficients arise in several
other situations in Physics and Engineering like in modeling
continuous sedimentation in clarifier thickener units used in waste
water treatment plants (See \cite{D1}, \cite{D2}, \cite{BKRT1}), in
traffic flow on highways with changing surface conditions (see
\cite{MOCHON}) and in ion etching used in the semiconductor industry
(see \cite{ROSS}). A detailed account of the above applications can be
found in \cite{SID1}. Consequently, equations of this type have been
studied extensively from both a theoritical as well as a numerical point
of view.  

In particular, for equations governing two-phase flow in porous media,
the numerical scheme that is commonly used by the petroleum engineers
is the upstream mobility flux scheme (see \cite{AzizSet,Sammon,BJ1,J1}). An
alternative finite difference (volume) method of the Godunov type
based on exact solutions of the Riemann problem was presented in
Adimurthi, Jaffre and Gowda (\cite{AJG1}). The above paper also
addressed the problem of prescribing correct entropy conditions at the
interface between rock types and showing the uniqueness of the entropy
solution. The 
solutions computed by the Godunov type scheme was shown to converge to
the entropy solution. The numerical flux developed in \cite{AJG1} is
similar to that used by Kaasschieter in \cite{KASS1} although written
in a more compact form. 

It is natural to ask whether the solution computed by the upstream
mobility flux scheme also converges to the entropy solution and
compare its numerical performance with other schemes like the one in
\cite{AJG1} and that of Towers (\cite{T1}, \cite{T2}). The goal of
this paper is to address these questions. 

In this paper, we will give an explicit representation of the upstream
mobility flux scheme for a medium consisting of two rock types and show
that the numerical flux is monotone. This will help us to obtain
estimates in $L^{\infty}$. We will then use a suitable modification
of the singular mapping technique to show that the approximate
solution converges to a weak solution of the continuous problem. . The
key point of this paper is to address whether this weak solution is a
entropy solution or not. We show by various numerical experiments that
the solutions computed by the upstream mobility flux scheme are not
consistent with the interface entropy condition of
\cite{AJG1}. Furthermore, we construct numerical experiments for the
case where only the absolute permeability changes. In this case, the solutions
given by the scheme do not converge to the entropy solution pointwise
and differ qualitatively from the entropy solution. The lack of
entropy consistency leads us to suggest that the upstream mobility
flux scheme may not be the correct numerical method to simulate two-phase flow in porous media with changing rock types and should be
replaced by the Godunov type scheme develop in \cite{AJG1}.

This paper is organised as follows, In section 2, we
describe the equations governing two phase flow in porous media with
heterogenities. In section 3, we summarise the mathematical theory
for single conservation laws with discontinuous flux developed
in \cite{AG1}, \cite{AJG1} and mention the well posedness
results. Section 4 is devoted to describing finite difference schemes
of the Godunov type as well as the upstream mobility flux scheme. The
convergence analysis for the upstream mobility flux scheme is carried
out in section 5. We give explicit representation formula for the
flux, show that it is monotone and use a variation of the singular
mapping technique to show convergence to a weak solution. The core of
this paper is in section 6 where we address the question of entropy
consistency. First, we consider numerical experiments for the case
where only the absolute permeability changes and discuss the entropy
behaviour of the solutions. Next, we construct examples of the case
where the relative permeability also changes and show that the scheme
is not consistent with the interface entropy condition of
\cite{AJG1}. We derive some conclusions from this paper in section 7.     

\section{Two-phase flow equations}
\setcounter{equation}{0}
Capillary-free two-phase incompressible flow is modeled by the following
scalar nonlinear hyperbolic equation
\[ \phi \frac{\p S}{\p t}+\frac{\p f}{\p x}=0 \]
where $\phi$ is the porosity of the rock, $S = S_1$ is the saturation of
phase 1 which lies in a bounded interval $[0,1]$.
The flux function $f$ is the Darcy velocity $\varphi_1$ of phase 1
and has the form
\[ f = \varphi_1 = \displaystyle{\frac{\lambda_1}{\lambda_1 + \lambda_2}
                  [ q + (g_1-g_2)\lambda_2 ]}.\]
Here $q = \varphi_1 + \varphi_2$ denotes the total Darcy velocity
where $\varphi_\ell, \ell=1,2$, denotes the Darcy velocity of phase $\ell$ with,
for the second phase,
\[ \varphi_2 = \displaystyle{\frac{\lambda_2}{\lambda_1 + \lambda_2}
 [ q + (g_2-g_1)\lambda_1 ]}. \]
Since the flow is assumed to be incompressible, the total Darcy
velocity $q$ is independent of the space variable $x$.

The quantities $\lambda_\ell, \, \ell=1,2$ may be called effective mobilities. 
These are products of the absolute permeability $K$ by the mobilities
$k_\ell$ :
\begin{displaymath}
\lambda_\ell = K k_\ell, \; \ell=1,2.
\end{displaymath}
The absolute permeability $K$ may depend on $x$ and the quantities $k_\ell$
and $\lambda_\ell$ are functions of $S$ which satisfy the following 
properties :
\[ \begin{array}{l}
k_1 \, \mbox{and} \, \lambda_1 \, \mbox{are increasing functions of} \, S ,\;
k_1(0) = \lambda_1(0) = 0 ,\\
k_2 \, \mbox{and} \, \lambda_2 \, \mbox{are decreasing functions of} \, S ,\;
k_2(1) = \lambda_2(1) = 0 .
\end{array} \]
We also shall assume that these functions are smooth functions of the
saturation $S$ and so is the flux function $f$.

The gravity constants $g_\ell, \, \ell=1,2$ of the phases are
\begin{displaymath}
g_\ell = g \rho_\ell \displaystyle{\frac{dx}{dz}}, \; \ell=1,2 \, ,
\end{displaymath}
with $g$ the acceleration due to gravity, $\rho_\ell$ the density of phase
$\ell$ and $z$ is the vertical position of the point of abscissa $x$.\\

Observe that with the above hypothesis, $f$ is a smooth (say
Lipschitz) function with at most one local maxi
with $f(0) = 0$ and $f(1) = q$ respectively.

\begin{figure}[htb]
\begin{center}
\begin{minipage}[c]{6cm}
\[
\left. \begin{array}{c}
\mbox{\bf Rock type I}\\ \\
\phi^{-},\; K^{-},\; k_1^{-},\; k_2^{-}\\
0 \leq S \leq 1\\
f \equiv f^{-}
\end{array} \right|
\begin{array}{c}
\mbox{ \bf Rock type II}\\ \\
\phi^{+},\; K^{+},\; k_1^{+},\; k_2^{+}\\
0 \leq S \leq 1\\
f \equiv f^{+}
\end{array}
\]
\end{minipage} $\qquad \qquad$
\begin{minipage}[c]{4cm}
\begin{picture}(100,100)(0,0)
\put(0,0){\includegraphics[height=4cm]{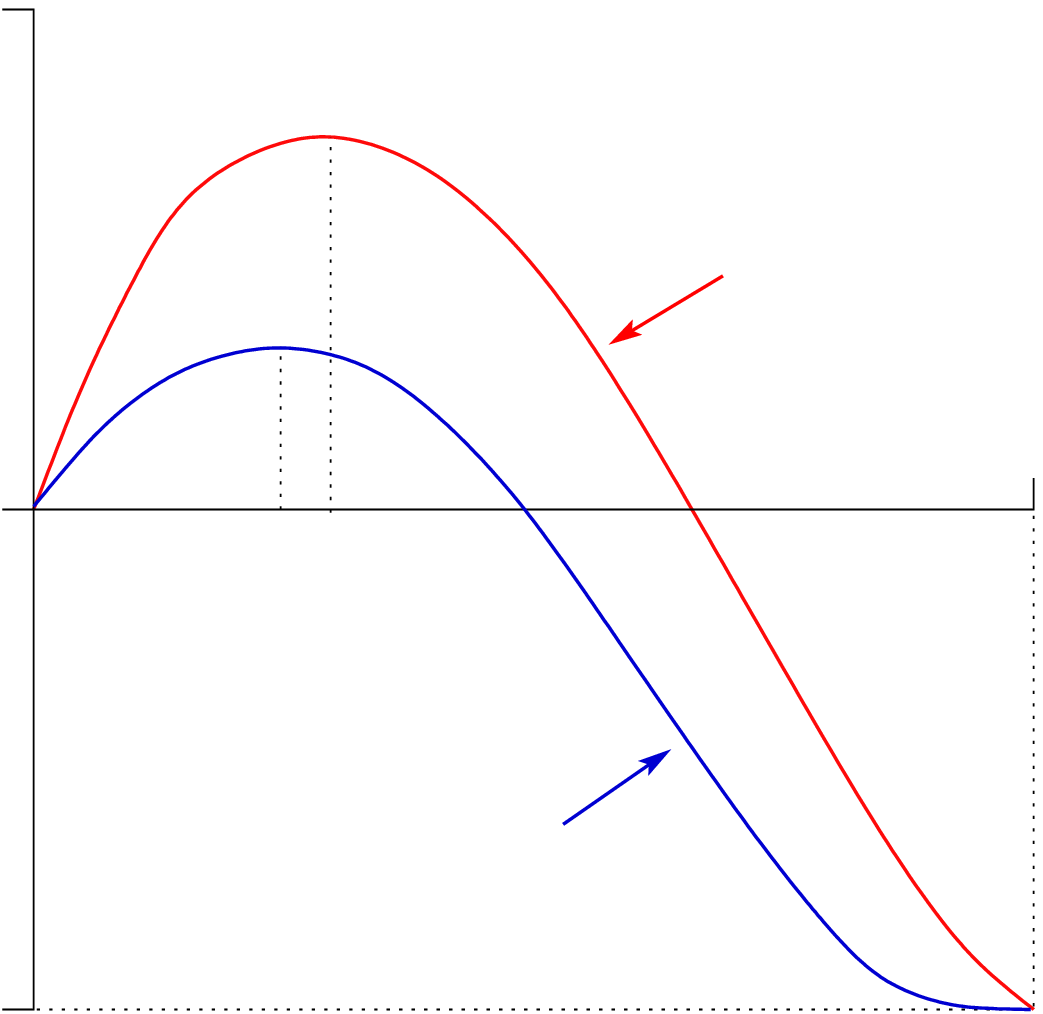}}
\put(55,20){\makebox(0,0){\textcolor{blue}{$f^-$}}}
\put(90,85){\makebox(0,0){\textcolor{red}{$f^+$}}}
\put(-5,0){\makebox(0,0){$q$}}
\put(-5,60){\makebox(0,0){0}}
\put(45,50){\makebox(0,0){\textcolor{red}{$\theta^+$}}}
\put(30,50){\makebox(0,0){\textcolor{blue}{$\theta^-$}}}
\end{picture}
\end{minipage}
\end{center}
\caption{Constants, mobilities and fluxes for two rock types}
\label{fig1}
\end{figure}

In many practical situations, the porous medium is heterogenous. For
example the medium may consist of two rock types separated
at the interface ($x=0$). In this case, the porosity, the absolute
permeability and the relative permeability change across the interface
and the flow is modeled by the following equations:
\begin{equation}
\label{eq4}
\begin{array}{l}
(H(x)\phi^{+} + (1-H(x))\phi^{-})S_{t} + (H(x)f^{+}(S) +
(1-H(x))f^{-}(S))_{x} = 0, \\
S(0,x) =  S_{0}(x) \end{array}
\end{equation}
where $H$ is the Heaviside function and the indices - and + refer to 
to the left and the right of the interface
respectively. The flux functions are given by the following,
\begin{equation}
\label{fluxequation}
f^{\pm} = \displaystyle{\frac{\lambda_1^{\pm}}{\lambda_1^{\pm} +
    \lambda_2^{\pm}}[ q + (g_2-g_1)\lambda_2^{\pm} ]}, \quad
\lambda_{i}^{\pm} = K^{\pm} k_{i}^{\pm}.  \end{equation}
See Figure \ref{fig1} for notations and the shapes of the flux functions when $q < 0, g_2 > g_1$. Note that the flux functions may also intersect as we will see in the numerical experiments of Section  \ref{numex}.

Equations (\ref{eq4}) with flux functions (\ref{fluxequation}) are a special case for the more general single
conservation law with flux function discontinuous in the space
variable considered in \cite{AJG1}. Flux functions $f^{-}$ and $f^{+}$ 
satisfy the following hypothesis,\\

\noi ${\bf H_1.} f^{-}, f^{+}$ are smooth (say Lipschitz) on $[0,1]$.\\
${\bf H_{2}.} f^{-}(0) = f^{+}(0) = 0$, $f^{-}(1)=f^{+}(1)=q$. \\
${\bf H_{3}.} f^{-},f^{+}$ have exactly one local maximum in $[0,1]$ with
$\theta^{-}= \mbox{ argmax }(f^-)$ and $\theta^{+} = \mbox{ argmax }(f^+)$. \\

Note that these are precisely the hypotheses for the fluxes assumed in \cite{AJG1}.

\section{The continuous problem}
As presented in the previous section, a change of rock types leads to a
single conservation law with a flux function discontinuous in the space
variable. An entropy theory has been developed for equations of the
form (\ref{eq4}) with the fluxes satisfying the hypothesis $H_{1},
H_{2}, H_{3}$. We summarize some of the results for the benefit of 
the reader.

Even in the case where the flux is continuous, it is well known that
solutions of equations of the above type develop discontinuities in
finite time even
for smooth initial data. Hence as in the continuous case, weak
solutions $S$ of (\ref{eq4}) are sought and defined as satisfying
\begin{multline}
\label{eq31}
\int_\R \int_{\R_{+}} S \varphi_t dx dt + \int_{\R} 
\int_{\R_{+}} (H (x) f^{+} (S) + (1- H (x) f^{-} (S)) \varphi_x dx dt \\
\hspace*{2cm} + \int_{\R} S_0 (x) \phi (0, x) dx = 0, \quad \forall \varphi \in
C_0^\infty (\overline{\R \times \R_+}).
\end{multline}
It is easy to check that $S$ is a weak solution of (\ref{eq4}) iff it satisfies in 
the weak sense, 
\[ \begin{array}{llll} 
S_t + (f^{-} (S))_x &=& 0, \quad x <0,\quad t >0, \\
S_t + (f^{+}(S))_x &=& 0, \quad x >0,\quad t >0, \\
S (0, x) &=& S_0 (x), \quad \forall \; x \in \R, \end{array} \]
and the following interface Rankine Hugoniot condition, 
\[
f^{+} (S^{+} (t)) = f^{-} (S^{-} (t)) \qquad \mbox{for almost all } t, \]
where
\[
S^+ (t) = \lim\limits_{x \rightarrow 0+} S (x, t), \hspace*{3cm}
S^- (t) = \lim\limits_{x \rightarrow 0-} S (x, t).
\]         

It is well known that weak solutions for a single conservation law are
not necessarily unique. Additional admissibility criteria termed as entropy
conditions need to be imposed for uniqueness. For equations of the form
(\ref{eq4}), it is natural to impose the standard Kruzkhov entropy
conditions away from the interface $x=0$. These can be stated in terms
of the entropy flux pairs which are defined as \\  
\noi {\bf Entropy pairs}:
$\{\varphi_{i}, \psi_{i}\}_{i=1,2}$ is said to be an  entropy pair for (\ref{eq4})
if $\varphi_{i}$ is convex in $[0, 1]$ and $\psi_1^\prime (\theta) =
\phi_1^\prime (\theta) f^{+ \prime} (\theta), \psi_2^\prime
(\theta) = \phi_{2}^{\prime} (\theta) f^{- \prime} (\theta) \mbox{ for } \theta \in [0,1].$

Let $S_0 \in L^\infty (\R)$ be the initial data with $0 \leq S_0
(x)\leq  1, \forall x \in \R,$ and let $S$ be a weak solution of (\ref{eq31}) with $0 \leq S(x,t)\leq  1, \forall (x,t) \in \R\times\R_+$.\\
\noi {\bf Interior entropy condition:} With $S_0$ and $S$ as above, $S$ is said to satisfy an
interior entropy condition  if for any entropy pairs $(\varphi_i,
\psi_i)_{i = 1,2}, S$ satisfies in the sense of distributions
\be
\label{32} 
\begin{array}{llll}
\dsp{\frac{\pa \varphi_1 (S)}{\pa t}} &+& \dsp{\frac{\pa \psi_1 (S)}{\pa x}}\leq 0, \quad \forall x >0, t >0. \\
\dsp{\frac{\pa \varphi_2 (S)}{\pa t}} &+& \dsp{\frac{\pa \psi_2 (S)}{\pa x}} \leq 0,
\quad \forall x < 0,t >0 \end{array} \ee

But for equation (\ref{eq4}), interior entropy conditions
like the one above are not sufficient to guarantee uniqueness and we need to
impose an additional entropy condition at the interface. The central
issue in the analysis of conservation laws with discontinuous flux is
the choice of this interface entropy condition. In \cite{AJG1}, the
following interface entropy condition was used.\\
\noi {\bf Interface entropy condition:} With $S_0$ and $S$ as above, assume
that $S^+ (t) = \lim\limits_{x \rightarrow 0^+} S (x,t)$ and $S^- (t) =
\lim\limits_{x \rightarrow 0^-} S (x, t)$  exist for almost all $t >0$
and define,
\[ \begin{array}{lcl}
L &=& \left\{t >0; \; S^- (t) \in (\theta^-, 1], S^+ (t) \in [0,
 \theta^+) \right\},  \\
U &=& \left\{t \in L; \; S^+ (t) =  S^-(t) = 1 \right\} \cup
\left\{t \in L; \; S^- (t) = S^+(t) = 0 \right\}.
\end{array} \]
Then $S$ is said to satisfy the interface entropy condition if
\begin{equation}
\label{33}
{\rm meas}\; \left\{L \setminus U \right\} = 0\,. \end{equation}
This  means that the characteristics must connect back to the
$x$-axis on at least one side of the jump in the flux i.e.,
undercompressive waves are not allowed. Undercompressive waves
i.e., ($f^{+ \prime}(S^{+}) > 0,\quad f^{- \prime}(S^{-}) < 0$) are
unrealistic physically as information is not taken from the initial
line. 

$S$ is defined to be an entropy solution of (\ref{eq4}) if it is a
weak solution and it satisfies both the interior as well as the interface
entropy condition. With this concept of entropy solution, the
following wellposedness result was obtained in \cite{AJG1}, under the
following hypotheses on the initial data:\\

\noi ${\bf IN_1.} \qquad  S_0 \mbox{ is such that } 0 \leq S_0 (x) \leq 1, 
\quad \forall \quad x \in  \R,$ \\
${\bf IN_2.} \qquad N (f^-, f^+, S_0) \leq  C < + \infty,$\\  

\noi where $N(f^-,f^+,S_0)$  is an estimator of the total variation of the
flux function evaluated at $S_0$. This estimator will be defined
precisely below in Section 4. 

We also need the following definition, \\
{\bf Regular solution.} $S$ is said to be a regular solution of
(\ref{eq4}) if the discontinuities of $S$ form a discrete set of
Lipschitz curves.

The well posedness result is given by 
\begin{theorem}
\label{theo2} 
Let $S_0$ satisfy hypotheses $IN_1,IN_2$ and $f^{-}, f^{+}$ satisfy 
hypotheses $H_{1}, H_{2}, H_{3}.$ Then there exists a weak solution $S \in L^\infty (\R \times \R+)$ of 
(\ref{eq4}) satisfying the following, 

\noi (1) \quad For almost all $t>0$  and $x \in \R, \quad S (x_-, t), 
S(x_+, t)$ exist. \\
\noi (2) \quad $S$ satisfies the interior entropy condition (\ref{32}).\\ 
\noi (3) \quad If $S$ is a regular solution, then $S$ satisfies the interface 
entropy condition (\ref{33}) and it is  unique.
\end{theorem} 

Uniqueness is proved by using a Kruzkhov type doubling of variables
argument. For details, see \cite{AJG1}. Existence was shown by showing
that a Godunov type finite difference scheme converges to a weak
solution and is consistent with the entropy conditions. 

We would also like to mention that more recently hypotheses on the fluxes have
been relaxed considerably to include fluxes of the concave-convex type
as in \cite{SID2, AGS2} ,  and with finitely many extrema as in
\cite{AGS2a}. Similarly schemes of the Enquist-Osher (EO) type have
been considered in \cite{AG1,SID2}. It is to be observed that
equations of the type (\ref{eq4}) are special cases of the more
general single conservation law with discontinuous coefficient of the
form, 
\[
u_{t} + (f(k(x),u))_{x} = 0, \qquad u(0,x) =  u_{0}(x) 
\]
The wellposedness and numerical methods for this problem are addressed
in a forthcoming paper \cite{AGS3}. We must also mention that an
entropy theory for equations of the above type (including a degenerate
parabolic term) has been developed by Karlsen, Risebro and Towers in
\cite{KRT1}. In \cite{T1,T2}, the author developed
staggered mesh algorithms of the Godunov and Enquist Osher type for
the multiplicative case i.e. $f(k,u) = k(x)f(u)$. The case with a
degenerate parabolic term was handled in \cite{KRT1} and the time
dependent case in \cite{CR1}. The entropy condition of \cite{KRT2}
agrees with that of \cite{AJG1} in many cases but differs in certain
cases. See Section 6 for a discussion of different entropy conditions
for equation (\ref{eq4}). For the rest of this paper, we will use the
entropy framework developed in \cite{AJG1}.

\section{Finite Difference Schemes}
In this section, we present finite difference schemes using for
numerical flux either the extended Godunov flux \cite{AJG1} or the
upstream mobility flux used in the petroleum industry for the
simulation of two-phase flow in porous media. 

Let $f$ be a Lipschitz continuous function. Then the Godunov flux
corresponding to $f$ is given by 
\begin{equation}
F_g(a, b) \; \;=\;\; \left\{ \begin{array}{rllll}
\min\limits_{\theta \in [a,b]} f (\theta) &{\rm if}& a< b, \\
\max\limits_{\theta \in [b,a]} f (\theta) &{\rm if}& a
\geq b, \end{array}
\right. \label{goddef} \end{equation}
The subscript $g$ stands here for Godunov in order to differentiate this flux from the upstream mobility flux that we will introduce later.
This flux was first proposed in \cite{G1} and is very popular in the
numerical analysis of conservation laws. It is based on exact
solutions of the Riemann problem. Let $F^{-}_{g}$ and $F^{+}_{g}$ be
the Godunov fluxes corresponding to the fluxes $f^{-}$ and $f^{+}$
respectively. 

In the case of two-phase flow, the flux functions $f^{-},f^{+}$
satisfy hypotheses ${\rm H_1, H_2, H_3}$ and the formula
(\ref{goddef}) can be simplified as follows,
\bea
F^{-}_{g}(a,b) &=& \min \{ f^{-}(\min (a,\theta_{-}), f^{-}(\max
(\theta_{-},b)\}, \nonumber \\
F^{+}_{g}(a,b) &=& \min \{ f^{+}(\min (a,\theta_{+}), f^{+}(\max
(\theta_{+},b)\}. \nonumber
\eea
These formulas were introduced in \cite{AJG1} and are simpler to implement than the general formula (\ref{goddef}).

Also, following \cite{AJG1}, they extend easily to define the
interface Godunov flux 
$\overline{F}_{g}$ based on the exact solution of the Riemann problem for
(\ref{eq4}): 
\bea
\label{interface}
\overline{F}_{g}(a,b) &=& \min \{ f^{-}(\min (a,\theta_{-}), f^{+}(\max
(\theta_{+},b)\}
\eea
We remark that the above interface flux coincides with the interface flux
obtained in \cite{KASS1} for which expression (\ref{interface})
represents a compact form, very easy to
use for computational purposes. It is also easy to check that the
interface flux $\ov{F}_{g}$ is Lipschitz is both variables,
nondecreasing in the first variable and nonincreasing in the second
variable. Note that the interface flux satisfies
\[ \overline{F}_{g}(0,0) = f^{-}(0,0) = f^{+}(0,0) = 0, \quad
\overline{F}_{g}(1,1) = f^{-}(1,1) = f^{+}(1,1) = q, \]
but is not consistent. 

Equipped with the definition of the numerical fluxes, we
proceed to describe the mesh. Let $h >0$ and define the space grid
points as follows: 
\[
x_{-1/2} = x_{1/2} = 0, \quad
x_{j + 1/2} = j\; h \quad {\rm for \;} j
\geq 0,\quad x_{j-1/2} =jh \quad
{\rm for}\; j \leq 0.
\]
We will also use the midpoints of the intervals:
\[
x_j = \left(\frac{2j-1}{2} \right) h
\quad {\rm for}\; j \geq 1, \quad
x_j = \left(\frac{2j+1}{2} \right) h
\quad {\rm for}\;  j \leq -1.
\]
For time discretization the time step is $\De t > 0$, and let
$t_n = n \De t, \; \la = \frac{\De t}{h}$.

For an initial data $S_0   \in L^\infty (\R)$  we define
$$
S^0_{j+1} = \frac{1}{h} \int_{x_{j+1/2}}^{x_{j+3/2}} \;
S_0 (x) dx \quad {\rm if} \quad j \geq 0, \quad
S^0_{j-1} = \frac{1}{h} \int_{x_{j-3/2}}^{x_{j-1/2}} S_0 (x)
dx \quad {\rm if}\quad j \leq 0.
$$

Now we can define  the  Godunov type finite difference scheme $\{S_j^n \} $
inductively  as follows:
\begin{equation}
\begin{array}{lcl}
\label{eq99}
S_1^{n+1} &=& S_1^n - \la (F^{+}_{g} (S_1^n, S_2^n)
- \overline{F}_{g} (S^n_{-1},S_1^n)),\\
S_j^{n+1} &=& S_j^n - \la (F^{+}_{g} (S_j^n, S^n_{j+1})
- F^{+}_{g} (S_{j-1}^n, S_j^n))\;\quad {\rm if} \; j > 1, \\
S_{-1}^{n+1} &=& S_{-1}^n - \la (\overline{F}_{g}
(S_{-1}^n, S_1^n ) - F^{-}_{g}(S_{-2}^n, S_{-1}^n)),  \\
S_j^{n+1} &=& S_j^n - \la (F^{-}_{g} (S_j^n,  S_{j+1}^n)
- F^{-}_{g} (S_{j-1}^n, S_j^n))\quad {\rm if}\quad j < -1.
\end{array}
\end{equation}
Observe that this is the standard Godunov scheme for  $j \neq \pm 1$, that
is, away from $x=0$, 

For $S_0 \in L^\infty (\R)$ and grid length $h$ and $\De t$ with
$\la = \frac{\De t}{h}$ fixed, define the piecewise constant function
$S_h \in L^\infty (\R \times \R_+)$ associated with $\{S_j^n\}$
calculated by the scheme (\ref{eq99}):
\begin{equation}
\label{defsh}
S_h (x, t) = S_j^n\,\; \mbox{ for } (x, t)
\in [x_{j-1/2}, x_{j+1/2} ) \times [n
\De t, \; (n+1) \De t), \quad j \neq 0.
\end{equation}

The above Godunov type scheme was analysed in \cite{AJG1}. For this analysis we need to introduce
\begin{eqnarray*}
N^{g}_h(f^{-},f^{+},S_0)&=&\sum\limits_{j
<-1}|F^{-}_{g}(S_j^0,S_{j+1}^0)-F^{-}_{g}(S_{j-1}^0,S_j^0)| \no \\
&&+
\sum\limits_{j > 1 } |F^{+}_{g}(S_j^0,S_{j+1}^0) -F^{+}_{g}(S_{j-1}^0,S_j^0)|\no \\
&&+ |\overline {F}_{g}(S_{-1}^0,S_1^0)-F_{g}^{-}(S_{-2}^0,S_{-1}^0)| \no \\
&&+|F^{+}_{g}(S_1^0,S_2^0)-\overline{F}_{g}(S_{-1}^0,S_1^0)|,
\end{eqnarray*}
\[
N_{g}(f^{-},f^{+},S_0) = \sup\limits_{h >0} N^{g}_h(f^{-},f^{+},S_0).
\]
It is easy to see that if $S_0 \in BV(\R)$, then $N_{g}(f^{-},f^{+},S_0)\leq
C|| S_0||_{BV}$, where $C$ is a constant depending only on the
Lipschitz constants of $f^{-}$ and $f^{+}$.

The following convergence theorem was proved, Let $M= \max Lip\{ f^{-},
  f^{+} \}$,
\begin{theorem}\label{theo1}
Assume that $\la, M$ satisfies the CFL condition $ \la M \leq 1.$
Let $ \; S_0 \in L^\infty(\R)  $ such that $0 \leq S_0 (x) \leq 1
\mbox{ for all }x \in \R$ and $N_{g}(f^{-},f^{+},S_0) <\infty.$ For $h
>0$, let $\la = \frac{\De t}{h}$ and $S_h$ be the corresponding
calculated solution given by {\rm (\ref{eq99}),
(\ref{defsh})}. Then there exists a subsequence $h_k \rightarrow
0$ such that $S_{h_k}$ converges almost everywhere to a weak solution $S$  of
{\rm (\ref{eq4})} satisfying the interior entropy condition. Suppose
the discontinuities of every limit function 
$S$ of $\{S_h\}$ form a discrete set of Lipschitz curves, then $S_h
\rightarrow S$ in $L^\infty_{loc} (\R_+, L^1_{loc} (\R))$ as $h
\rightarrow 0$, and $S$ satisfies the interface  entropy condition.
\end{theorem}

The convergence of the scheme was proved by using the singular mapping
technique which we will also use in section 5 albeit with
modifications. The limit solution obtained was shown to be consistent
with the interior as well as the interface entropy condition. The key
point in the proof of consistency with the interface entropy condition
was the use of a contradiction argument using s test function. The
reader is referred to \cite{AJG1} for details. We will use similar
ideas in the next section. Some numerical experiments involving this
Godunov type scheme are shown in section 6.

As mentioned earlier, staggered mesh schemes were proposed in
\cite{T1}, \cite{T2} and \cite{KRT1} for general single conservation
laws with discontinuous flux. In the simplified case of a single
discontinuity in the flux, the staggered mesh scheme of the Godunov
type can also be
written in the form (\ref{eq99}) by replacing the interface Godunov flux
$\ov{F}_{g}$ with the averaged interface flux $\ov{F}_{\tau}(a,b)$ which is
the Godunov flux corresponding to the function $\tau
=1/2(f^{-}+f^{+})$. This finite difference scheme is analyzed
in \cite{T1}  and is shown to converge for a large class of
fluxes. Numerical experiments comparing this scheme with a Godunov
type scheme was reported in \cite{SID2}. We will also compare this
scheme in the numerical experiments in Section 6.

The main objective of this paper is to analyse the upstream mobility
flux. It is an adhoc flux for two phase flow in porous media, invented
by petroleum engineers from simple physical considerations, and it
corresponds to an approximate solution of the Riemann problem. The
standard upstream mobility flux for $f^{-}$ is given by the following
formula:
\begin{equation}
\label{stanup-}
\begin{array}{l}
F^{-}(a,b)= \displaystyle{
         \frac{\lambda_1^{-*}}{\lambda_1^{-*} + \lambda_2^{-*}}
         [ q + (g_1-g_2)\lambda_2^{-*} ]} ,\\
\lambda^{-*}_\ell=
\left\{ \begin{array}{ll} \lambda^{-}_\ell(a) & \mbox{if }
         q+(g_\ell-g_i)\lambda_i^{-*} >0, \;i=1,2, i\neq \ell,\\[3mm]
        \lambda^{-}_\ell(b) & \mbox{if }
        q+(g_\ell-g_i)\lambda_i^{-*} \leq0, \; i=1,2, i\neq \ell,
\end{array}
\right. \ell=1,2,
\end{array}
\end{equation} 

Similarly, the standard upstream mobility flux corresponding to $f^{+}$
can be defined by the following formula,
\begin{equation}
\label{stanup+}
\begin{array}{l}
F^{+}(a,b)= \displaystyle{
         \frac{\lambda_1^{+*}}{\lambda_1^{+*} + \lambda_2^{+*}}
         [ q + (g_1-g_2)\lambda_2^{+*} ]} ,\\
\lambda^{+*}_\ell=
\left\{ \begin{array}{ll} \lambda^{+}_\ell(a) & \mbox{if }
         q+(g_\ell-g_i)\lambda_i^{+*} >0, \;i=1,2, i\neq \ell,\\[3mm]
        \lambda^{+}_\ell(b) & \mbox{if }
        q+(g_\ell-g_i)\lambda_i^{+*} \leq0, \; i=1,2, i\neq \ell,
\end{array}
\right. \ell=1,2,
\end{array}
\end{equation} 
These formulas just say that the mobility $\lambda_\ell$ must be
calculated using the value of the saturation which is upstream with
respect to the flow of the phase $\ell$ since the sign of the
quantity $q+(g_1-g_2)\lambda_2^{*}$ determines the direction of the
flow of phase 1 and the sign of $q+(g_2-g_1)\lambda_1^{*}$
determines that of phase 2.

Note that the above formulae are implicit and have been made explicit
in \cite{BJ1}. The flux is shown to be Lipschitz, monotone and
consistent in the same reference. 

As for the Godunov scheme, an
interface upstream mobility scheme needs to be defined to take into
account the changing rock types. Formulas (\ref{stanup-}),(\ref{stanup+})
can be easily extended to obtain the interface flux $\ov{F}(a,b)$:
\begin{equation}
\label{stanupi}
\begin{array}{l} 
\ov{F}(a,b)= \displaystyle{
         \frac{\lambda_1^{*}}{\lambda_1^{*} + \lambda_2^{*}}
         [ q + (g_1-g_2)\lambda_2^{*} ]} ,\\
\lambda^{*}_\ell=
\left\{ \begin{array}{ll} \lambda^{-}_\ell(a) & \mbox{if }
         q+(g_\ell-g_i)\lambda_i^{*} >0, \;i=1,2, i\neq \ell,\\[3mm]
        \lambda^{+}_\ell(b) & \mbox{if }
        q+(g_\ell-g_i)\lambda_i^{*} \leq0, \; i=1,2, i\neq \ell,
\end{array}
\right. \ell=1,2,
\end{array}
\end{equation}
This interface flux preserves the idea of calculating the flux using
the phase mobilities which are upstream with respect to the flow of the
corresponding phases. 

Now we define the upstream mobility flux scheme for a medium
with changing rock types as follows,
\begin{equation}
\begin{array}{lcl}
\label{eq100}
S_1^{n+1} &=& S_1^n - \la (F^{+} (S_1^n, S_2^n)
- \overline{F} (S^n_{-1},S_1^n)),\\
S_j^{n+1} &=& S_j^n - \la (F^{+} (S_j^n, S^n_{j+1})
- F^{+} (S_{j-1}^n, S_j^n))\;\quad {\rm if} \; j > 1, \\
S_{-1}^{n+1} &=& S_{-1}^n - \la (\overline{F}
(S_{-1}^n, S_1^n ) - F^{-}(S_{-2}^n, S_{-1}^n)),  \\
S_j^{n+1} &=& S_j^n - \la (F^{-} (S_j^n,  S_{j+1}^n)
- F^{-} (S_{j-1}^n, S_j^n))\quad {\rm if}\quad j < -1.
\end{array}
\end{equation}

For $S_0 \in L^\infty (\R)$ and grid length $h$ and $\De t$ with
$\la = \frac{\De t}{h}$ fixed, define the function
$S_h \in L^\infty (\R \times \R_+)$ associated with $\{S_j^n\}$
calculated by the scheme (\ref{eq100}):
\begin{equation}
\label{defSh}
S_h (x, t) = S_j^n\,\; \mbox{ for } (x, t)
\in [x_{j-1/2}, x_{j+1/2} ) \times [n
\De t, \; (n+1) \De t), \quad j \neq 0.
\end{equation}

We will analyse the scheme (\ref{eq100}) in the next section. As for
the Godunov case we will need a $BV$ type norm which we define as  
\begin{eqnarray*}
N_h(f^{-},f^{+},S_0)&=&\sum\limits_{j
<-1}|F^{-}(S_j^0,S_{j+1}^0)-F^{-}(S_{j-1}^0,S_j^0)| +
\sum\limits_{j > 1 } |F^{+}(S_j^0,S_{j+1}^0) -F^{+}(S_{j-1}^0,S_j^0)|\no \\
&&+ |\overline {F}(S_{-1}^0,S_1^0)-F^{-}(S_{-2}^0,
S_{-1}^0)|+|F^{+}(S_1^0,S_2^0)-\overline{F}(S_{-1}^0,S_1^0)|,
\end{eqnarray*}
\[
N(f^{-},f^{+},S_0) = \sup\limits_{h >0} N_h(f^{-},f^{+},S_0).
\]

\section{Convergence Analysis}
In this section, we show that the solutions defined by
(\ref{eq100}),(\ref{defSh}) converge to a weak solution of (\ref{eq4})
along a subsequence as $h \rightarrow 0$. We closely follow the
analysis of \cite{AJG1} and will refer to the above paper for
details. We first observe that the formulae (\ref{stanup-}),
(\ref{stanup+}), (\ref{stanupi}) are implicit. The first step is to make
them explicit. For the interior fluxes, this has been done in
\cite{BJ1}. We will give an explicit representation of the interface
flux. Depending on the ordering of the gravity constants, we have to
distinguish the following two cases.\\

\noi {\bf Case 1}: $g_{1} \leq g_{2}$\\
\noi Following \cite{BJ1}, we define the  auxillary quantities for
calculating the explicit fluxes
\[ \begin{array}{lcllcl}
\theta_{1}&=& q +(g_{1}-g_{2})\lambda^{-}_{2}(a), & \quad
\delta_{1}&=& q +(g_{1}-g_{2})\lambda_{2}^{*}, \nonumber \\
\theta_{2}&=& q +(g_{2}-g_{1})\lambda^{+}_{1}(b),  & \quad 
\delta_{2}&=& q +(g_{2}-g_{1})\lambda_{1}^{*}. \nonumber
\end{array} \]
\noi Clearly we have $\theta_{1} \leq \theta_{2}$ and
$\delta_{1} \leq \delta_{2}$. We have the following lemma for the
explicit formulae of the fluxes,
\begin{lemma}
\label{lemma-E1}We can have only the following three cases:
\[ \begin{array}{lrcl}
1. & 0 \leq \theta_1 = \delta_{1} \leq \theta_{2} &\Leftrightarrow&
\lambda_{1}^{\ast} = \lambda^{-}_{1}(a),
\lambda_{2}^{\ast}=\lambda^{-}_{2}(a)\\
2. &\theta_1 = \delta_{1} \leq 0 \leq \theta_2 = \delta_{2} &\Leftrightarrow&
\lambda_{1}^{\ast} = \lambda^{+}_{1}(b),
\lambda_{2}^{\ast}=\lambda^{-}_{2}(a)\\
3. &\theta_{1} \leq \theta_2 = \delta_{2} \leq 0 &\Leftrightarrow&
\lambda_{1}^{\ast} = \lambda^{+}_{1}(b),
\lambda_{2}^{\ast}=\lambda^{+}_{2}(b)
\end{array} \]
\end{lemma}
The proof is simple and similar to the case of the interior flux
$F$. Details can be found in \cite{SID1}. This lemma says that just
calculating $\theta_1$ and $\theta_2$ is sufficient to determine the
upstream side of the flow of each phase.

The other case works in the same way.\\

\noi {\bf Case 2}: $g_{2} \leq g_{1}$\\
We define the auxillary quantities
\[
\theta_{1} = q +(g_{1}-g_{2})\lambda^{+}_{2}(b), \quad
\theta_{2} = q +(g_{2}-g_{1})\lambda^{-}_{1}(a)
\]
and $\delta_1$ and $\delta_2$ are as in case 1. Now we have
 $\theta_{1} \geq \theta_{2}$ and
$\delta_{1} \geq \delta_{2}$ and the  following lemma.
\begin{lemma}
\label{lemmaE2} We can have only the following three cases
\[ \begin{array}{lrcl}
1. &\theta_{1} \geq \theta_{2}= \delta_2 \geq 0 &\Leftrightarrow&
\lambda_{1}^{\ast} = 
\lambda^-_1(a),\quad \lambda_{2}^{\ast}=\lambda^{-}_{2}(a)\\
2. &\theta_{1}= \delta_1\geq 0 \geq \theta_{2}=\delta_2 &\Leftrightarrow&
 \lambda_{1}^{\ast} = 
\lambda^{-}_{1}(a),\quad \lambda_{2}^{\ast}=\lambda^{+}_{2}(b)\\
3. & 0 \geq \theta_{1}=\delta_1 \geq \theta_{2} &\Leftrightarrow&
 \lambda_{1}^{\ast} = 
\lambda^{+}_{1}(b),\quad \lambda_{2}^{\ast}=\lambda^{+}_{2}(b)
\end{array} \]
\end{lemma}

Again calculating $\theta_1$ and $\theta_2$ gives the direction of the
flow of each phase and the way to calculate the upstream mobilities.
This is easy to implement in a code. 

Our goal is to show that the sequence of approximate saturations converges to a
weak solution of (\ref{eq4}). We begin by stating some of the
properties of the interface flux.
\begin{lemma}
\label{lemma-P}
The interface flux $\ov{F}$ as defined in (\ref{stanupi}) is Lipschitz
in both its arguments, non decreasing in the first and nonincreasing
in the second argument. Furthermore the following also holds
\[ \ov{F}(0,0) = f^{-}(0) = f^{+}(0) = 0, \quad
\ov{F}(1,1) = f^{-}(1) = f^{+}(1) = q. \]
\end{lemma}
{\it Proof}: The proof that the flux is Lipschitz is
similar to that for the interior fluxes in \cite{BJ1} and we
omit the details. Also the evaluation of $\ov{F}(0,0)$ and
$\ov{F}(1,1)$ is easy to check. So let us have a quick pass at the
monotonicity properties. We consider for example case 1 i.e. $g_1 \leq
g_2$.

If $0 \leq \theta_1 = \delta_{1} \leq \theta_{2}$ we have
\[ \begin{array}{l}
\lambda_{1}^{\ast} = \lambda^{-}_{1}(a),\quad
\lambda_{2}^{\ast} = \lambda^{-}_{2}(a), \quad 0 \leq 
\delta_1 = q+(g_1-g_2)\lambda^{-}_{2}(a) \leq 
\delta_2 = q+(g_2-g_1)\lambda^{-}_{1}(a),\\[0.3cm]
\ov{F}(a,b) = \dsp{\frac{\lambda^{-}_{1}(a)}
{\lambda^{-}_{1}(a)+\lambda^{-}_{2}(a)}}\,\delta_1, 
\quad \dsp{\frac{\p \ov{F}}{\p a}(a,b)} = \dsp{\frac{
(\lambda_1^-)^{\prime}(a) \lambda^{-}_{2}(a) \delta_1 + 
\lambda_1^-(a) (\lambda^{-}_{2})^{\prime}(a)(-\delta_2)}
{(\lambda^-_1(a)+\lambda^{-}_{2}(a))^2}}.
\end{array} \]
Since $\lambda^{-}_{1}, \lambda^{-}_{2}$ are both positive functions,
$\lambda^{-}_{1}$ is nondecreasing and $\lambda^{-}_{2}$ is
nonincreasing, and $0 \leq \delta_1 \leq \delta_2$, we conclude that
$\dsp{\frac{\p \ov{F}}{\p a}(a,b)} \geq 0$ and $\ov{F}(a,b)$ is
nondecreasing with respect to $a$. Obviously $\ov{F}(a,b)$ does not
depend on $b$ so it is nonincreasing with respect to $b$.

If $\theta_1 = \delta_{1} \leq 0 \leq \theta_2 = \delta_{2}$ we have
\[ \begin{array}{l}
\lambda_{1}^{\ast} = \lambda^{+}_{1}(b), \quad
\lambda_{2}^{\ast} = \lambda^{-}_{2}(a), 
\quad \delta_1 = q+(g_1-g_2)\lambda^{-}_{2}(a) \leq 0 \leq
\delta_2 = q+(g_2-g_1)\lambda^{+}_{1}(b), \\[0.3cm]
\ov{F}(a,b) = \dsp{\frac{\lambda^{+}_{1}(b)}
{\lambda^{+}_{1}(b)+\lambda^{-}_{2}(a)}}\, \delta_1, 
\; \dsp{\frac{\p \ov{F}}{\p a}(a,b)} = \dsp{\frac{
\lambda_1^+(b) (\lambda^{-}_{2})^{\prime}(a)(- \delta_2)}
{(\lambda^{+}_{1}(b)+\lambda^{-}_{2}(a))^2}}, \;
\dsp{\frac{\p \ov{F}}{\p b}(a,b)} = \dsp{\frac{
 \lambda^{-}_{2}(a)(\lambda^{+}_{1})^{\prime}(b) \delta_1}
{(\lambda^{+}_{1}(b)+\lambda^{-}_{2}(a))^2}}.
\end{array} \]
Again it is easy to check that $\dsp{\frac{\p \ov{F}}{\p a}(a,b)} \geq 0$ and 
$\dsp{\frac{\p \ov{F}}{\p b}(a,b)} \leq 0$. Therefore  $\ov{F}(a,b)$
is nondecreasing with respect to $a$ and nonincreasing with respect
to $b$. 

If $\theta_1 = \delta_{1} \leq \theta_2 = \delta_{2} \leq 0$ we have
\[ \begin{array}{l}
\lambda_{1}^{\ast} = \lambda^{+}_{1}(b), \quad
\lambda_{2}^{\ast} = \lambda^{+}_{2}(b), 
\quad \delta_1 = q+(g_1-g_2)\lambda^{+}_{2}(b) \leq
\delta_2 = q+(g_2-g_1)\lambda^{+}_{1}(b) \leq 0 , \\[0.3cm]
\ov{F}(a,b) = \dsp{\frac{\lambda^{+}_{1}(b)}
{\lambda^{+}_{1}(b)+\lambda^{+}_{2}(b)}}\, \delta_1, 
\quad \dsp{\frac{\p \ov{F}}{\p b}(a,b)} = \dsp{\frac{
(\lambda_1^+)^{\prime}(b) \lambda^{+}_{2}(b) \delta_1 + 
\lambda_1^+(b) (\lambda^{+}_{2})^{\prime}(b)(-\delta_2)}
{(\lambda^+_1(b)+\lambda^{+}_{2}(b))^2}}.
\end{array} \]
Again it is easy to check that $\dsp{\frac{\p \ov{F}}{\p a}(a,b)} = 0$ and 
$\dsp{\frac{\p \ov{F}}{\p b}(a,b)} \leq 0$ and $\ov{F}(a,b)$
is nondecreasing with respect to $a$ and nonincreasing with respect
to $b$. \cqfd\\

Similar statements for the interior fluxes can be found in
\cite{BJ1}. 

Next, we state the CFL condition for stability of the
scheme as the following,
\begin{equation}
\label{CFL}
\begin{array}{ll}
\lambda M \leq 1, \\
M = \max \left\{ \right.
&\max\limits_{|j| > 1,n} 
\{\dsp{\frac{\pa F^{n}_{j+1/2}}{\pa a}(S_{j},S_{j+1}) - 
\frac{\pa F^{n}_{j-1/2}}{\pa b}(S_{j-1},S_j)}\}, \\
&\left.\dsp{\frac{\pa F^{n}_{3/2}}{\pa a}(S_1,S_2) - \frac{\pa \ov{F}^{n}}{\pa b}(S_{-1},S_1),
\frac{\pa \ov{F}^{n}}{\pa a}(S_{-1},S_1) - 
\frac{\pa F^{n}_{-3/2}}{\pa b}(S_{-2},S_{-1})}
\right\}. \end{array} \end{equation}
 
This type of a condition was explicitly written out in \cite{Sammon} for the
case of one rock type and a slight modification of it gives the result
in our case.
Then we can prove in a straightforward manner 
\begin{lemma}
\label{lemma-M}
Under the CFL condition (\ref{CFL}), the upstream mobility scheme
defined by (\ref{eq100}) is monotone.
\end{lemma}

We remark that scheme (\ref{eq100}) is in conservative
form, is monotone, but it is not consistent because of the interface
flux (some examples are discussed in section 6). Hence, the
classical theory (see \cite{CM, KV, GR1}) does not apply and we have to adopt
the analysis presented in \cite{AJG1}. The monotonicity of the scheme
leads to the following discrete $L^{1}$ contractivity result:

\begin{lemma}
\label{discrete}
Let $S_0, \in  L^\infty (\R,\; [0,1])$ be the initial data, and
let $\{S_j^n\}$ be the corresponding solution
calculated by the upstream mobility flux scheme
{\rm(\ref{eq100})}. then,
\begin{equation}
\label{discrete1}
\sum\limits_{j \neq 0} |S_j^{n+1} - S_j^{n}| \leq
\sum\limits_{j \neq 0} |S_j^n - S_j^{n-1}|.
\end{equation}
\end{lemma}
{\it Proof}: As the scheme (\ref{eq100}) is monotone and
conservative, this estimate follows by applying the Crandall-Tartar
lemma (see \cite{GR1}). \cqfd \\

The next step is to obtain estimates in $L^{\infty}$ for the
approximate solutions. For Monotone, Consistent and Conservative
schemes , such estimates follow from a discrete maximum principle. But
for the scheme (\ref{eq100}), the lack of consistency implies that the
discrete maximim principle is no longer true. Instead, as in
\cite{AJG1}, we can use the consistency of the interface flux at the
points $0$ and $1$ to obtain that $[0,1]$ is an invariant region for
the scheme and obtain the following lemma,
\begin{lemma}
\label{maxbounds}
Let $S_0 \in L^\infty (\R,[0,1])$ be the initial data, and let
$\{S_j^n\}$ be the corresponding solution calculated by the
finite volume scheme {\rm (\ref{eq100})}. The following holds,
\begin{equation}
\label{0S1}
0 \leq S_j^n \leq 1 \;\ \forall j, n.
\end{equation}
\end{lemma}
{\it Proof.} Since $0 \leq S_0 \leq 1$, hence for all $j, \; 0
\leq S_j^0
\leq 1.$ By induction, assume that (\ref{0S1}) holds
for $n$. Then from Lemma \ref{lemma-P} we have
$$
\begin{array}{llll}
0 &=& H_{-1} (0, 0, 0) \leq H_{-1} (S_{j-1}^n,
S^n_{j}, S_{j+1}^n) = S_1^{n+1}
\leq H_{-1} (1,1, 1) = 1 \; {\rm if} \; j \leq -2, \\
0 &=& H_1 (0, 0, 0) \leq H_1 (S_{j-1}^n, S_j^n,
S_{j+1}^n) = S_j^{n+1} \leq H_1 (1,1, 1) = 1
\quad {\rm if } \; j \geq 2, \\
0 &=& H_{-2} (0, 0, 0) \leq H_{-2} (S_{-2}^n,
S_{-1}^n, S_1^n) = S_{-1}^{n+1}
\leq H_{-2} (1,1, 1) = 1,\\
0 &=& H_2 (0, 0, 0) \leq H_2 (S_{-1}^n ,
S_1^n, S_2^n) = S_1^{n+1} \leq H_2
(1, 1, 1) = 1.
\end{array}
$$
This proves (\ref{0S1}). \cqfd \\

As pointed out earlier, the key difficulty in the convergence analysis
is to obtain $BV$ type estimates on the approximations. For a
monotone, consistent and conservative scheme, such estimates following
from the discrete $L^{1}$ contractivity and the translation
invariance (see \cite{GR1}). But in this case, the scheme is not
consistent and we cannot expect the approximate solutions to be
uniformly bounded in $BV$. Rather, the difficulty is circumvented by
using the singular mapping technique first introduced by Temple in
\cite{Temple1} and adapted for schemes for single conservation laws by
Towers in \cite{T1}. The singular mapping was also adopted in the
convergence proof in \cite{AJG1}. More recently, several modifications
of the singular mapping have been suggested in \cite{SID2},
\cite{AGS1} etc.

The central idea in using the singular mapping technique is to
estimate the total variation of the approximate solutions under the
singular mapping by the variation of the fluxes in neighboring cells
and use the discrete $L^{1}$ contractivity. This method works well for
upwind schemes like Godunov and Enquist Osher but it does not work for
other types of numerical fluxes like the Lax-Friedricks flux. The
same is true for the upstream mobility flux and we have to adapt the
technique to work in this case. We do so by using the idea of chain
estimates like in \cite{AGS2} . We start by defining the
singular mappings. ,we use the following standard notation $a \in \R,\mbox{then}\quad a_{+} = \max\{a,0\},a_{-} = \min\{a,0\},a=a_{+}+a_{-},|a|=a_{+}-a_{-}.$\\
\noi The singular mappings are given by,
\begin{equation}
\label{sm-1}
\psi_{1}(\theta) = \int_{\alpha}^{\theta}|f^{- \prime}(\xi)|d\xi, \quad
\psi_{2}(\theta) = \int_{\alpha}^{\theta}|f^{+ \prime}(\xi)|d\xi
\end{equation}
where $\alpha \in [0,1]$ is some number. Note that we use the
same form of singular mappings as in \cite{AJG1} expect that there are centerered
at an arbitrary point as in \cite{SID2}. Now we are in a position to
define the transformed schemes for the discrete values of the
solution. We define them as 
\[
z_{j}^{n} = \left\{ \begin{array}{lcl}
\psi_{1}(S_{j}^{n}) & {\rm if } & j \leq -1 \\
\psi_{1}(S_{-1}^{n}) & {\rm if} & j \geq -1 
\end{array} \right. ,
w_{j}^{n} =  \left\{ \begin{array}{lcl}
\psi_{2}(S_{1}^{n}) & {\rm if } & j \leq 1 \nonumber \\
\psi_{2}(S_{j}^{n}) & {\rm if } &  j \geq 1 
\end{array} \right. .
\]

Like in \cite{SID2}, we define two sets of transformed variables which
enables us to simplify the proof to some extent as compared to
\cite{AJG1}. Our goal is to estimate the variation of the transformed
scheme at each time level. For simplicity, let us suppres the
subscript $n$ as we are dealing with the same time level. Then
\bea
TV(z_{j}) = \sum\limits_{j \neq 0}|z_{j}-z_{j+1}| = 2 \sum\limits_{j
  \neq 0}(z_{j}-z_{j+1})_{+} \nonumber 
\eea
In \cite{AJG1},\cite{SID2}, this
variation was controlled individually in each cell by the flux
variation across the neighboring cells. For details see lemma (5.4) in
\cite{SID2}. But such an estimate relied on the upwind nature of the
Godunov flux and is not necessarily true for the upstream mobility
flux as the upstream mobility flux gives different answers from the
Godunov and Enquist-Osher fluxes when the phases are flowing in
different directions. Rather, nonlocal variation estimates hold in
this case as will be explained below. For this, we observe from the
definition of the singular mapping (\ref{sm-1}) that
$(z_{j}-z_{j+1})_{+} > 0$ if and only if $S_{j} > S_{j+1}$. Same
observation also applies to $w_{j}$'s. We use
this ordering of the neighboring cell values to define the
following,\\
\noi Define  $\un{J}= \{ j \leq -2 \}$ and we define some subsets of this set as
follows,\\
\noi {\bf Definition}: $\un{I} \subset \un{J}= \{i \in \un{J} : S_{i} < S_{i+1}\}$ is the set of admissble
  indices.\\
\noi Note that this implies that for each $i \in \un{I}$, there exists a
unique $k(i)$ such that the following holds,\\
\noi 1. $S_{i} \leq S_{i-1} \leq \ldots \leq S_{i-k(i)}$\\
\noi 2. $S_{i-k(i)-1} < S_{i-k(i)}$ \\
\noi We denote the following,\\
\noi 1.$k(i) = 0$ if $S_{i-1} < S_{i}$ and \\
\noi 2.$k(i) = \infty$ if $\forall j < i,\quad S_{j} \geq S_{j+1}$\\
\noi So there can be atmost one $i \in \un{I}$ such that $k(i) = \infty$. Let $i_{0}$ be such that $i_{0} = \min\limits_{\un{I}} i$. Note that $i_{0}$ is not necessarily equal to $-2$.  Now we
can define a chain as $\un{J}_{i} = \{ j: k(i) \leq j \leq i \}$. With the
above definitions, it easy to check that $\un{J} =
\cup_{i \in \un{I}} \un{J}_{i}$.\\   

Similarly denote, $\ov{J}= \{ j \geq 1 \}$ and we define some subsets of this set as
follows,\\
\noi {\bf Definition}: $\ov{I} \subset \ov{J}= \{i \in \ov{J} : S_{i} < S_{i+1}\}$ is the set of admissble
  indices.\\
\noi Note that this implies that for each $i \in \ov{I}$, there exists a
unique $k(i)$ such that the following holds,\\
\noi 1. $S_{i} \geq S_{i+1} \geq \ldots \geq S_{i+k(i)}$\\
\noi 2. $S_{i+k(i)+1} > S_{i+k(i)}$ \\
\noi We denote the following,\\
\noi 1.$k(i) = 0$ if $S_{i+1} > S_{i}$ and \\
\noi 2.$k(i) = \infty$ if $\forall j > i,\quad S_{j} \geq S_{j+1}$\\
\noi So there can be atmost one $i \in \ov{I}$ such that $k(i) = \infty$. We denote the minimum value in $\ov{I}$ to be $i^{0}$. Now we
can define a chain as $\ov{J}_{i} = \{ j: i \leq j \leq k(i) \}$. With the
above definitions, it easy to check that $\ov{J} =
\cup
_{i \in \ov{I}} \ov{J}_{i}$.\\
\noi Equipped the definitions above, we are in a position to state the
main lemma of this section in the following,
\begin{lemma}
\label{fluxvar}
$\forall i \in \un{I}$, the following estimates hold, if $i < -2$ and
$k(i) \neq \infty$ then,
\bea
\label{fluxvar1}
\sum\limits_{j \in \un{J}_{i}}(z_{j} - z_{j+1})_{+} &\leq& \sum\limits_{j
  \in \un{J}_{i}}\{ |F^{-}(S_{j},S_{j+1}) - F^{-}(S_{j-1},S_{j})|
  \nonumber \\ 
&+&
|F^{-}(S_{j+1},S_{j+2}) - F^{-}(S_{j},S_{j+1})| \} 
\eea
for $i_{0}$ as defined above, and such that $i_{0} = -2$, we have the following estimate,
\bea
\label{fluxvar2}
\sum\limits_{j \in \un{J}_{i_{0}}}(z_{j} - z_{j+1})_{+} &\leq& 2M
\eea
In case, $i$ happens to be the only index such that $k(i) = \infty$,
  then 
\bea
\label{fluxvar3}
\sum\limits_{j \in \un{J}_{i}}(z_{j} - z_{j+1})_{+} &\leq& 2M
\eea
Similarly, if $i \in \ov{I}$ and $k(i) \neq -\infty$, then the following estimate holds
\bea
\label{fluxvar4}
\sum\limits_{j \in \ov{J}_{i}}(w_{j} - w_{j+1})_{+} &\leq& \sum\limits_{j
  \in \ov{J}_{i}}\{ |F^{+}(S_{j},S_{j+1}) - F^{+}(S_{j-1},S_{j})|
  \nonumber \\ 
&+&
|F^{+}(S_{j+1},S_{j+2}) - F^{+}(S_{j},S_{j+1})| \} 
\eea 
for $i^{0}$ as defined above, and such that $i^{0} = 1$, we have the following estimate,
\bea
\label{fluxvar5}
\sum\limits_{j \in \ov{J}_{i^{0}}}(z_{j} - z_{j+1})_{+} &\leq& 2M
\eea
And in case $i$ is such that $k(i)= \infty$, then we have that
\bea
\label{fluxvar6}
\sum\limits_{j \in \ov{J}_{i}}(w_{j} - w_{j+1})_{+} &\leq& 2M
\eea
\end{lemma}
{\it Proof}: We will only provide proofs for the estimates
(\ref{fluxvar1}) and (\ref{fluxvar3}). The other inequalities follow in
  the same manner. We have to consider three separate cases to check
  the estimate namely,\\
\noi {\bf Case 1}: $0 \leq S_{k(i)} \leq \theta^{-}$.\\
\noi In this case, it follows from the $L^{\infty}$ bounds and the
definitions that $0 \leq S_{i} \leq \ldots \leq S_{k(i)} \leq
\theta^{-}$. Hence, one can check that 
\bea
\label{eq513}
\sum\limits_{j \in \un{J}_{i}}(z_{j} - z_{j+1})_{+} &=&
f^{-}(S_{k(i)}) - f^{-}(S_{i}) 
\eea
\noi From the definition of $\un{I}$, we get that $S_{i} \leq S_{i+1}$
and $S_{k(i) + 1} \leq S_{k(i)}$. Therefore, using the monotonicity
and consistency of the interior upstream mobility flux scheme,we get
that
\bea
\label{eq514}
F^{-}(S_{i},S_{i+1}) \leq F^{-}(S_{i},S_{i}) = f^{-}(S_{i}) 
\eea
and
\bea
\label{eq515}
F^{-}(S_{k(i)},S_{k(i)+1}) \geq F^{-}(S_{k(i)},S_{k(i)}) = f^{-}(S_{k(i)}) 
\eea
Therefore by combining the above estimates, we get that 
\bea
f^{-}(S_{k(i)}) - f^{-}(S_{i}) &=&f^{-}(S_{k(i)}) -
F^{-}(S_{k(i)+1},S_{k(i)+2}| \nonumber \\
&+& F^{-}(S_{k(i)+1},S_{k(i)+2}) - \ldots \nonumber \\
&+& \ldots -F^{-}(S_{i-1},S_{i})\nonumber \\
&+& F^{-}(S_{i-1},S_{i}) - f^{-}(S_{i}) \nonumber
\eea 
Now by using (\ref{eq514}) and (\ref{eq515}), we get that ,
\bea
\sum\limits_{j \in \un{J}_{i}}(z_{j} - z_{j+1})_{+} &\leq& \sum\limits_{j
  \in \un{J}_{i}} |F^{-}(S_{j},S_{j+1}) - F^{-}(S_{j-1},S_{j})|
\nonumber
\eea
thus proving (\ref{fluxvar1}). Next we consider,\\
\noi {\bf Case 2}: $\theta^{-} \leq S_{i} \leq 1$.\\
\noi In this case, it follows from the $L^{\infty}$ bounds and the
definitions that $\theta^{-} \leq S_{i} \leq \ldots \leq S_{k(i)} \leq
1$. Hence, one can check that 
\bea
\label{eq516}
\sum\limits_{j \in \un{J}_{i}}(z_{j} - z_{j+1})_{+} &=&
f^{-}(S_{i}) - f^{-}(S_{k(i)}) 
\eea
\noi From the definition of $\un{I}$, we get that $S_{i} \leq S_{i-1}$
and $S_{k(i) - 1} \leq S_{k(i)}$. Therefore, using the monotoniticity
and consistency of the interior upstream mobility flux scheme,we obtain
\bea
\label{eq517}
F^{-}(S_{i-1},S_{i}) \geq F^{-}(S_{i},S_{i}) = f^{-}(S_{i}) 
\eea
and
\bea
\label{eq518}
F^{-}(S_{k(i)-1},S_{k(i)}) \leq F^{-}(S_{k(i)},S_{k(i)}) = f^{-}(S_{k(i)}).
\eea
Therefore by combining the above estimates, we have
\bea
f^{-}(S_{k(i)}) - f^{-}(S_{i}) &=&f^{-}(S_{k(i)}) -
F^{-}(S_{k(i)+1},S_{k(i)+2})| \nonumber \\
&+& F^{-}(S_{k(i)+1},S_{k(i)+2}) - \ldots \nonumber \\
&+& \ldots -F^{-}(S_{i-1},S_{i})\nonumber \\
&+& F^{-}(S_{i-1},S_{i}) - f^{-}(S_{i}) \nonumber
\eea 
Now by using (\ref{eq517}) and (\ref{eq518}), we get 
\bea
\sum\limits_{j \in \un{J}_{i}}(z_{j} - z_{j+1})_{+} &\leq& \sum\limits_{j
  \in \un{J}_{i}}\{ |F^{-}(S_{j},S_{j+1}) - F^{-}(S_{j-1},S_{j})|
\nonumber
\eea
thus proving (\ref{fluxvar1}).\\
\noi {\bf Case 3}:$S_{i} \leq \theta^{-} \leq S_{k(i)}$\\
\noi In this case, $\exists \quad l(i) \in \un{J}_{i}$ such that $S_{i} \leq
\ldots \leq S_{l(i)} \leq \theta^{-} \leq S_{l(i) - 1} \leq \ldots
\leq S_{k(i)}$. We get that
\bea
\label{eq519}
\sum\limits_{j \in \un{J}_{i}}(z_{j} - z_{j+1})_{+} &=&
f^{-}(\theta^{-})-f^{-}(S_{k(i)})+f^{-}(\theta^{-}) - f^{-}(S_{i}) 
\eea
Again by the monotonicity and consistency of the
interior fluxes, we have the following estimate,
\bea
\label{eq520}
F^{-}(S_{l(i)-1},S_{l(i)}) \geq F^{-}(\theta^{-},\theta^{-}) =
f^{-}(\theta^{-})
\eea
Therefore, from the estimates (\ref{eq519}), (\ref{eq520}), (\ref{eq514})
and (\ref{eq518}), we have 
\bea
f^{-}(\theta^{-}) - f^{-}(S_{i}) &=&f^{-}(\theta^{-}) -
F^{-}(S_{l(i)},S_{l(i)+1} \nonumber \\
&+& F^{-}(S_{l(i)},S_{l(i)+1}) - \ldots \nonumber \\
&+& \ldots -F^{-}(S_{i-1},S_{i})\nonumber \\
&+& F^{-}(S_{i-1},S_{i}) - f^{-}(S_{i}) \nonumber
\eea 
and
\bea
f^{-}(\theta^{-}) - f^{-}(S_{k(i)}) &=&f^{-}(\theta^{-}) -
F^{-}(S_{l(i)-2},S_{l(i)-1}) \nonumber \\
&+& F^{-}(S_{l(i)-2},S_{l(i)-1}) - \ldots \nonumber \\
&+& \ldots -F^{-}(S_{k(i)},S_{k(i)+1})\nonumber \\
&+& F^{-}(S_{k(i)},S_{k(i)+1}) - f^{-}(S_{k(i)}) \nonumber
\eea 
Combining the above 2 estimates, we get the desired inequality and
prove (\ref{fluxvar1}) in all the 3 cases. Next, we prove
(\ref{fluxvar3}). In case of $i$ being the unique element of $\un{I}$
such that $k(i) = \infty$. It is easy to from the $L^{\infty}$ and
Lipschitz bounds that 
\bea
\label{eq521}
\sum\limits_{-\infty}^{i} (z_{j} - z_{j+1})_{+} &\leq& f(\theta^{-}) - f^{-}(u_{i})\nonumber \\
&+& f^{-}(\theta^{-}) - f^{-}(1)\nonumber \\
&\leq & M(|\theta^{-}-1| + \theta^{-}) = M \nonumber
\eea
\noi Thus we have the estimate (\ref{fluxvar3}). The other estimates can be proved similarly. \cqfd \\

We use the above inequalities to show the following variation bound on the singular mapping,
\begin{lemma}
\label{lemma-Var}
The transformed sequences are of bounded total variation and the following estimate holds,
\bea
\max\{ TV(z_{j}^{n}),TV(w_{j}^{n}) \} &\leq& \frac{4}{\lambda} N(f^{-},f^{+},S_{0}) + 2M 
\eea
\end{lemma}

\noi{\it Proof}: We provide a proof for the sequence $\{ z_{j}^{n} \}$, the other bound follows in a similar way.We  
\bea
TV(z^{n}_{j}) &=& 2(\sum\limits_{j \in \Z}(z_{j}^{n} - z_{j+1}^{n})_{+}) \nonumber \\
&=&  \sum\limits_{j \in \un{J}}(z_{j}^{n} - z_{j+1}^{n})_{+} \nonumber \\
&=&  \sum\limits_{i \in \un{I}}\sum\limits_{j \in \un{J_{i}}}(z_{j}^{n} - z_{j+1}^{n})_{+} + 2M\nonumber \\
\eea
By adding the chain inequalities (\ref{fluxvar1}), \ref{fluxvar2} and \ref{fluxvar3}), we get that
\bea
\sum(z_{j}^{n} -z_{j+1}^{n})_{+} &=& 2(\sum\limits_{j \leq -2}|G(u_{j}^{n},u_{j+1}^{n}) -G(u_{j-1}^{n},u_{j}^{n})| + |G(u_{-2}^{n},u_{-1}^{n} - \overline{F}(u_{-1}^{n},u_{1}^{j})| \nonumber \\
&+& |\overline{F}(u_{-1}^{n},u_{1}^{j})- F(u_{1}^{n},u_{2}^{n})|+\sum\limits_{j \geq 2}|F(u_{j}^{n},u_{j+1}^{n}) -F(u_{j-1}^{n},u_{j}^{n})| + 2M \nonumber \\
&=& \frac{2}{\lambda}\sum\limits_{j \neq 0}|u_{j}^{n +1}-u_{j}^{n}| + 4M
\eea
Now by using the discrete $L^{1}$ contractivity, we get the desired estimate. In a similar way, we can get the total variation bound for $\{ w_{j}^{n} \}$ and complete the proof of the lemma. \cqfd \\

In order to show convergence of solutions generated by the scheme, we need to define the following piecewise constant functions,Let $z^{h},w ^{h}$ be defined as $z^{h}(x,t)= z_{j}^{n}, w^{h}(x,t)=w^{n}_{j},\quad \forall \quad (x,t) \in I_{j}^{n}$. We translate the bounds on the discrete values in terms of the above functions in the following lemmas which we state without proof (for a proof see \cite{SID2}).
\begin{lemma}
\label{lemma-Func}
With the functions defined as above and $\forall t \in \R_{+}$,we have,
\bea
\max\{TV(z^{h}),TV(w^{h})\} &\leq& \frac{4}{\lambda} N_{h}(f^{-},f^{+},u_{0}) + 4M
\eea
\end{lemma}
\begin{lemma}
\label{lemma-Func1}
Let $S_0 \in L^\infty (\R, [0, 1])$ such that $N(f^{-},f^{+},S_0) <
 \infty$ be initial data and let $S_h$  be the co
rresponding
solutions obtained by the scheme (\ref{eq4}),then
\bea
0 \leq S_h (x,t) \leq 1 \qquad \forall\; \; (x, t) \in \R \times \R_+ 
\eea
\bea
\int\limits_{\R} |S_h (x,t) - S_h (x, \tau)| dx &\leq& N_h(f^{-},f^{+},S_0)(2 \De t + |t-\tau|) 
\eea
\end{lemma}
We are in a position to state our main convergence theorem. The key step was to prove the total variation bounds on the singular mappings and the fact that the singular mapping is monotone and hence invertible. We have that

\begin{theorem}
\label{convergence}
Assume that the CFL condition is satisfied and the initial data $S_{0}$ satisfies the hypothesis $IN_{1}$ and $IN_{2}$, Let $S_{h}$ be approximate solutions defined above,then there exists a subsequence (still denoted by h) such that $S_{h}$ converge almost everywhere to a weak solution S of (\ref{eq4}).In fact $S_{h} \rightarrow S$ in $L^{\infty}_{{\rm loc}}(\R_{+},L^{1}_{{\rm loc}}(\R))$ as h goes to 0. Furthermore, the limit solution satisfies the interior entropy condition (\ref{32}).
\end{theorem}

\noi{\it Proof}:This is the main convergence theorem for the upstream mobility flux scheme (\ref{eq100}). The proof follows by the classical arguments of the Lax-Wendroff theorem (See \cite{GR1}) and the modifications introduced in 
\cite{AJG1}. We omit the details and refer to the above quoted paper for them.
\par For any fixed $t > 0$,we have the BV bounds from the Lemma (\ref{lemma-Func}) and by using the standard Rellich compactness theorem that upto subsequences (still denoted by $h$),we get that $z_{h}(.,t),w_{h}(.,t)$ converge in $L^{1}$ and for almost all x to $z(.,t)$ and $w(.,t)$ respectively.
\par For fixed t and for almost all $x < 0$, then from the convergence we have that $\psi(S_{h}) \rightarrow z$. Note that $\psi$ is monotone for $x < 0$ and we get that $S_{h}(x,t) \rightarrow \psi^{-1}(z(x,t))=S(x,t)$. Thus for all $x < 0$,we get that $S_{h}(.,t)$ converges to $S(.,t)$ for almost all $x < 0$.Similarly,for $x > 0$,we define $S(x,t) = \psi_{2}^{-1}(z_{h}(x,t)$ and get the a.e convergence.Now we use the standard density argument along with the time continuity estimate (\ref{lemma-Func1}),to get that 
\bea
S_{h} \rightarrow  S \quad \in L^{\infty}_{{\rm loc}}((0,T),L^{1}_{{\rm loc}}(\R))
\eea

\par Once we have the above convergence, we can use the standard
arguments of the Lax-Wendroff type in order to show that $S_{h}$
converges to a weak solution of (\ref{eq4}).The proof follows exactly
as in \cite{SID2} and we refer the reader to this reference for
details. Similarly, the consistency of the scheme with the interior
entropy condition is shown by using the numerical entropy fluxes of
Crandall-Majda (see \cite{GR1}). Check \cite{AJG1} for the details.  \cqfd

\section{Entropy consistency of the scheme} 
\label{numex}

In the previous sections, we have shown that the upstream mobility
flux scheme is well defined for two-phase flow in an heterogenous medium
with two rock types and that it generates solutions which converge to weak
solutions of the conservation law (\ref{eq4}) and satisfy  interior
Kruzkhov type entropy condition. But, for the solutions of the scheme
to be admissible, we have to show that they also satisfy the interface
entropy condition (\ref{33}). As pointed out earlier, we remark
that this entropy condition is pointwise and essentially amounts to
the exclusion of undercompressive waves at the interface ($x=0$). In
(\cite{AJG1}), it was shown that the Godunov type scheme (\ref{eq99})
satifies the interface entropy condition by means of a contradiction
argument. In this section, we investigate the question of whether the
limit solution generated by the upstream mobility flux scheme also
satisfies this interface entropy condition or not.  

It will be shown by means of various numerical experiments that the
limit solution generated by scheme (\ref{eq100}) which uses the
upstream mobility flux need not satisfy the
interface entropy condition and  counterexamples are given. In other
cases that we report, it is far from clear whether the interface
entropy condition is actually satisfied. In fact, the numerical evidence
suggests that this condition is not satisfied in the pointwise sense
as required in the entropy theory of (\cite{AJG1}) on account of
certain boundary layer phenomena at the interface although the entropy
condition may be satisfied in a weaker integral sense. 

It is worth mentioning that there are several entropy theories for
equations of the type (\ref{eq4}) like those presented in \cite{AJG1}
and in \cite{KRT2}. It will be shown that in some cases, the solutions
generated by the scheme (\ref{eq100}) satisfy the entropy conditions
of \cite{KRT2} and in some other cases, the conditions of
\cite{AJG1}. 

We will now present the five numerical experiments illustrating five
different situations.\\ 

\noi{\bf Experiment 1} 

\begin{figure}[H]
\begin{center}
\epsfysize=6cm
\epsffile{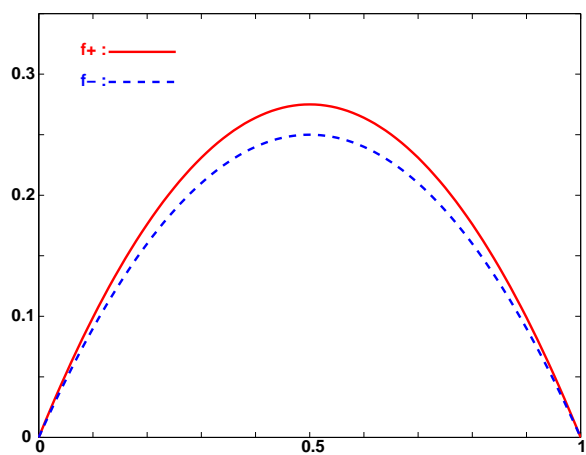}
\caption{Flux functions in experiment 1}
\end{center}
\label{figureE1}
\end{figure}

In this example, we consider the flux functions given by,
\[
\begin{array}{llllllllll}
\lambda_{1}^{+}(S) &=& 1.1S && \lambda_{2}^{+}(S)&=&1.1(1-S) \\
\lambda_{1}^{-}(S) &=& S && \lambda_{2}^{-}(S) &=& 1-S\\
g_{1} &=& 2 && g_{2}&=& 1 \\
\phi &=& 1 && q &=& 0
\end{array}
\]
As is clear from the above, we are considering that the porosity and
the relative permeabilities dont change across the rock types and the
absolute permeabilities only change with $K^{+}= 1.1$ and
$K^{-}=1$. The shape of the corresponding fluxes $f^{-}$ and $f^{+}$
are shown in Fig. \ref{figureE1}.
 
We consider for initial data,
$S_{0}(x)= \left\{ \begin{array}{lcl}
 0.65 &{\rm if}& x < 0 \\
 0.35 &{\rm if}& x > 0.
\end{array} \right. $
In this case, since the flux functions do not intersect in the interior
of the interval (0.1), the entropy solution in this case coincides
for the entropy theories of \cite{AJG1} and \cite{KRT2} and consists
of a rarefaction fan joining $0.65$ and $0.5$ on the left and a steady
discontinuity at the interface with $0.5$ as the left trace and $0.35$
as the right trace. Note that the entropy solution does not admit
underconmpressive waves at the interface as $f^{-\prime}(0.5) \equiv
0$.   

\begin{figure}[H]
\epsfysize=4.5cm 
\epsffile{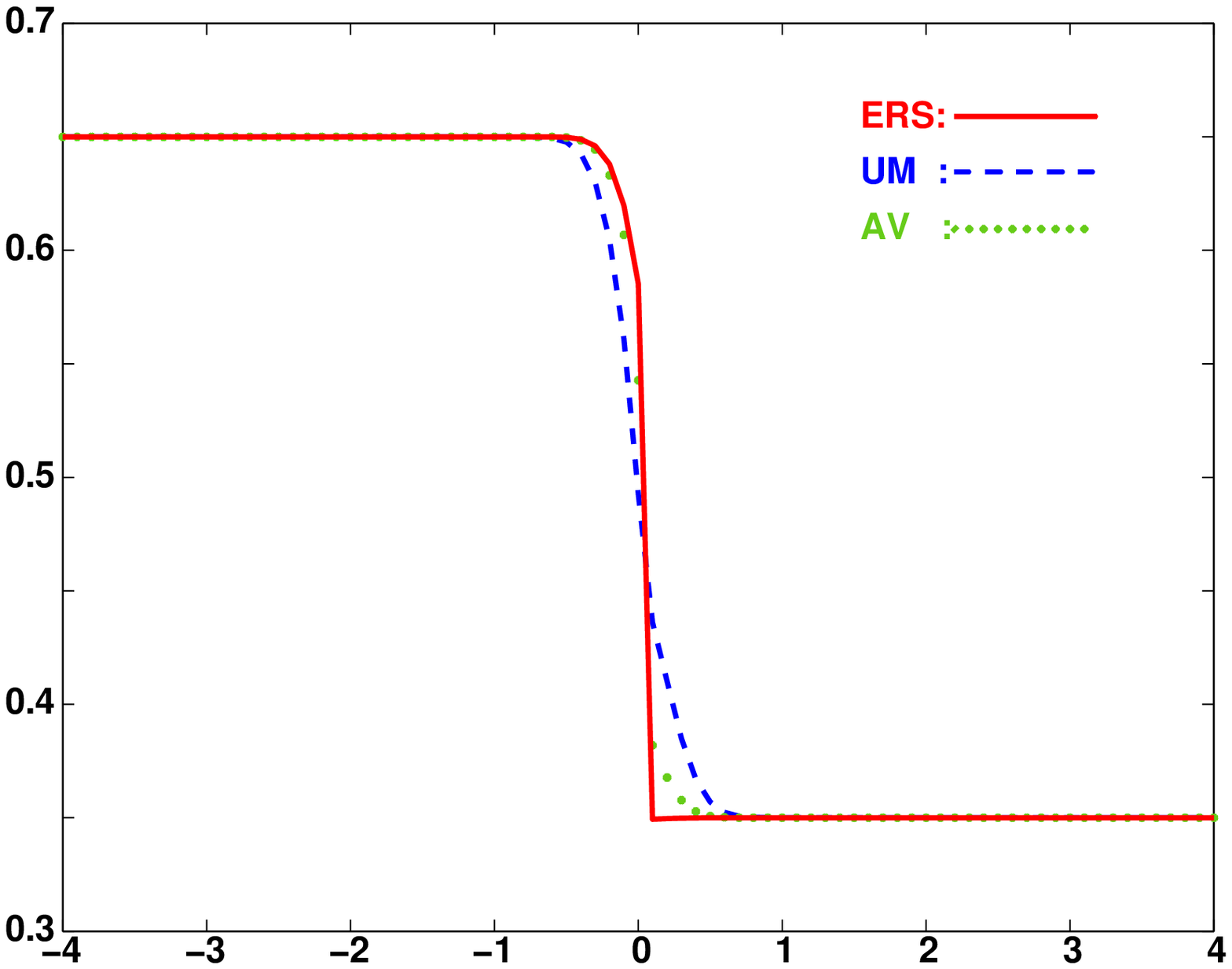} \hspace*{0.3cm} \epsfysize=4.5cm 
\epsffile{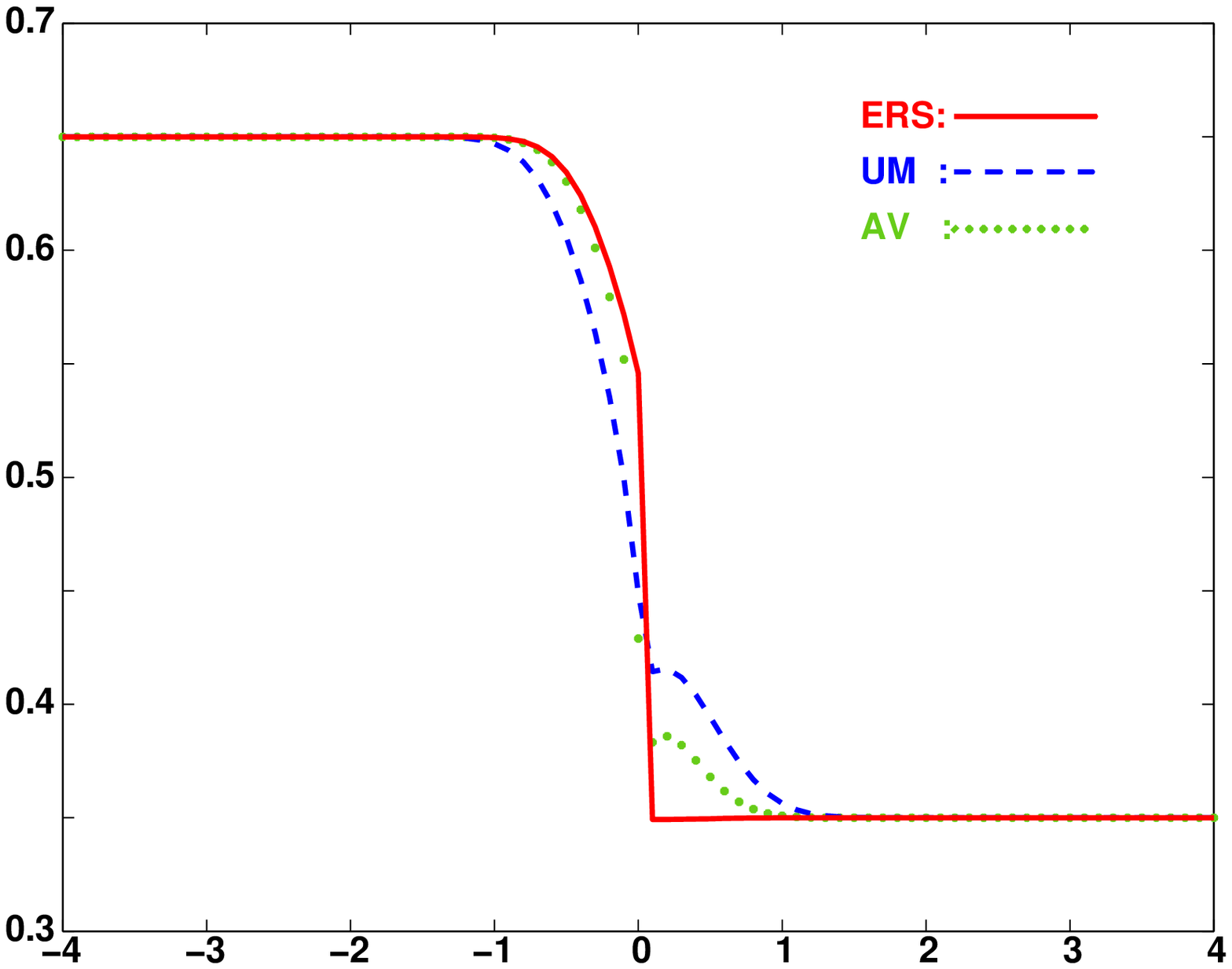}
\caption{Solutions in experiment 1 with h = 0.1 at times t=0.5 and t=1.5}
\label{figureE2}
\end{figure}

\begin{figure}[H]
\epsfysize=4.5cm 
\epsffile{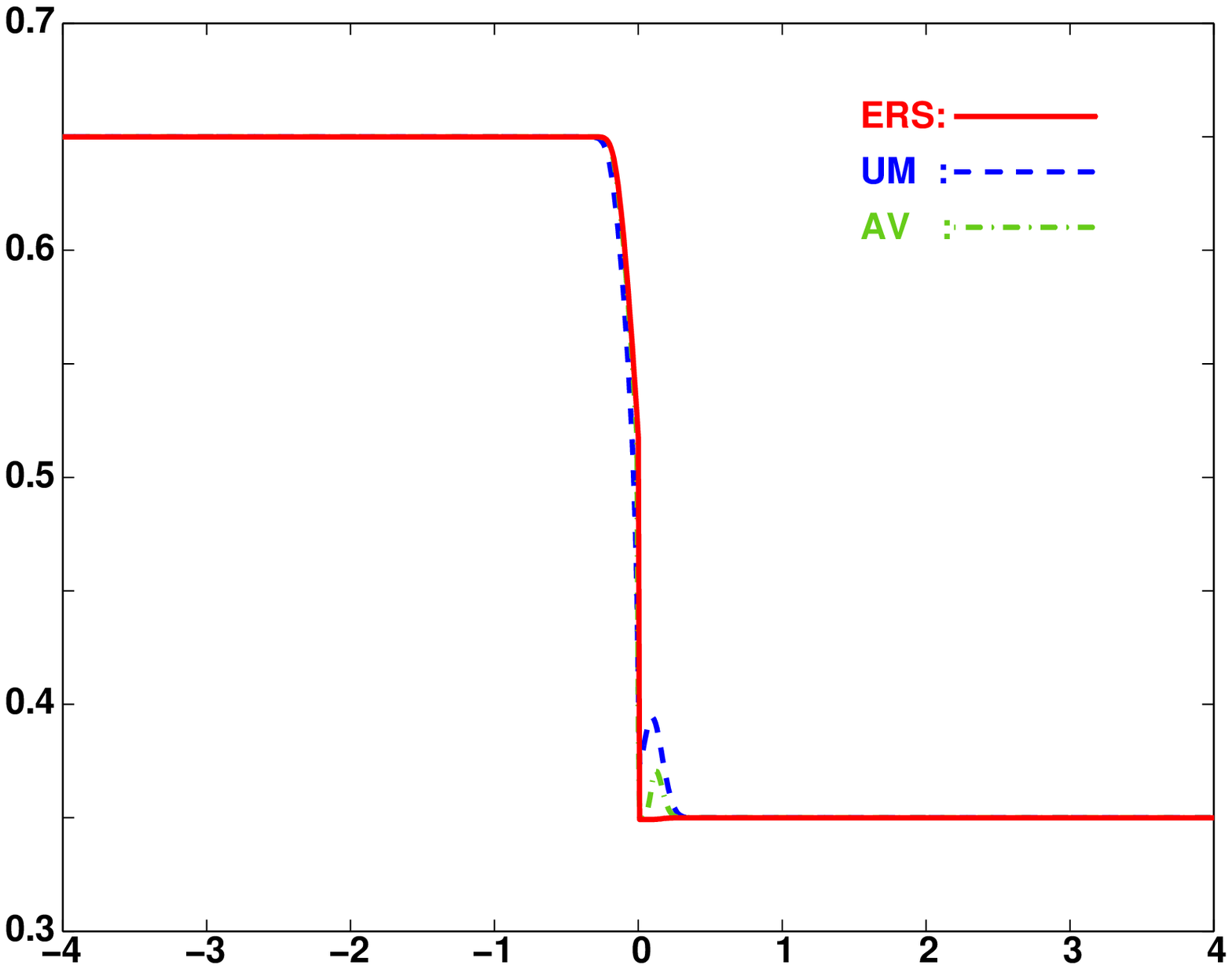} \hspace*{0.3cm} \epsfysize=4.5cm 
\epsffile{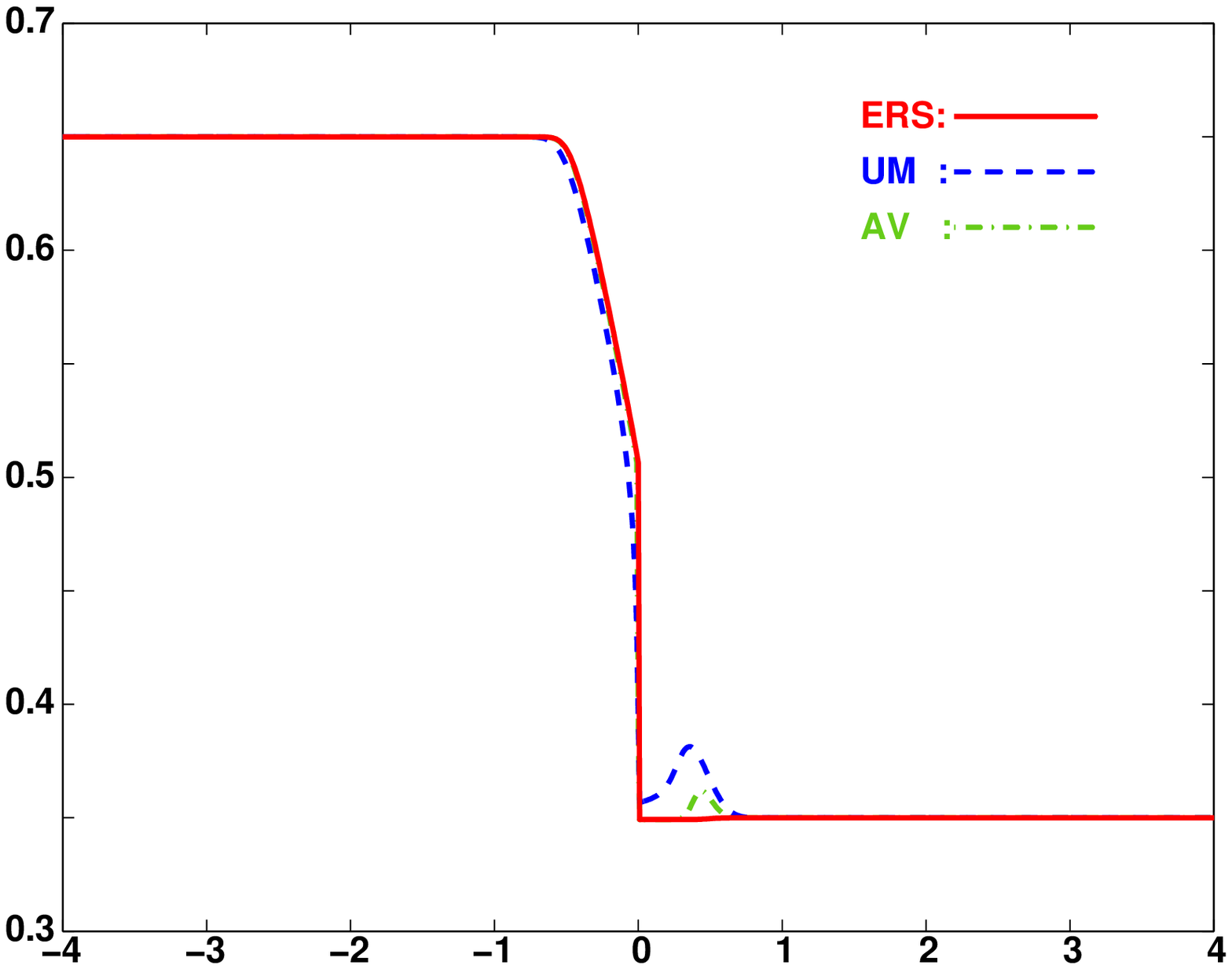}
\caption{Solutions in experiment 1 with h = 0.01 at t=0.5 and t=1.5}
\label{figureE3}
\end{figure}

We present the solutions obtained by the Godunov type finite
difference scheme which we term as the Exact Riemann Solver (ERS) and
the Upstream Mobility scheme which we term as UM. Also, we compute
solutions with the staggered mesh version of the Godunov scheme
developed by Towers in \cite{T1} and \cite{T2}. We term this scheme as
AV. The computed solutions are shown at two different times and for
two different mesh sizes in Fig. \ref{figureE2} and
Fig. \ref{figureE3}. 

In Figure \ref{figureE2}, we show the numerical results obtained
with the mesh size $h=0.1$ and the CFL constant
$\lambda=\frac{1}{8}$. As expected at this rather large mesh size, the resolution
is not high although the ERS is already giving very good results
with the interface discontinuity being resolved perfectly. On the
other hand, both UM and AV do not resolve the interface discontinuity
well. In fact, as seen in Fig. \ref{figureE2}, the left trace of the
solution as computed by UM is approximately $0.4$ which is less than
the expected trace of $0.5$. This is indicative of the development of
a boundary layer at the interface $x=0$. Another anomaly is the
existence of a travelling wave in both UM and AV that is clearly
unphysical as the solution is constant $0.35$ in $x > 0$. The
amplitude of this spurious wave is higher for UM than for AV. Both these
phenomena indicate that we cannot prove that the limit solution
generated by UM and AV are consistent with the interface entropy
condition (\ref{33}) in a pointwise sense. 

In order to confirm the above proposition, we reduce the mesh size to
$h=0.01$ and show the solution in Fig. \ref{figureE3}. Again, we
see that ERS resolves both the rarefaction and the interface
discontinuity very well with little numerical diffusion whereas both
UM and AV do not match the solution. Even with a very small mesh
size, the left trace at the interface of the solution computed with
UM is around $0.4$ and is well below the required value of
$0.5$. Also, the spurious travelling wave seen before is still present
although its magnitude has decreased. As stated earlier, the existence
of both a
boundary layer and a travelling wave forces us to believe that the
solution computed with UM is not consistent with the interface
entropy condition. The same holds true for the
solutions computed by AV. \\

The first numerical experiment that we have presented represents the
simplest type of discontinuity at the interface involving only a
change in the absolute permeability across the interface. Even in this
simple situation, the UM flux scheme does not perform as
well as ERS and the limit solution obtained by it doesnot seem to
satisfy the interface entropy condition of \cite{AJG1}. Hence, more
interesting and complicated behaviour is 
expected when we consider changes in relative permeabilities across
the interface. As will be shown in the coming numerical
experiments, the limit solution computed by UM will converge to the
entropy solution of \cite{AJG1} in some cases and the entropy solution
of \cite{KRT2} in some other cases. 

We start with an example where the solution given by UM seems to
converge to the entropy solution of \cite{AJG1} in the following
numerical experiment,\\ 
\noi {\bf Experiment 2}

\begin{figure}[H]
\begin{center}
\epsfysize=6cm
\epsffile{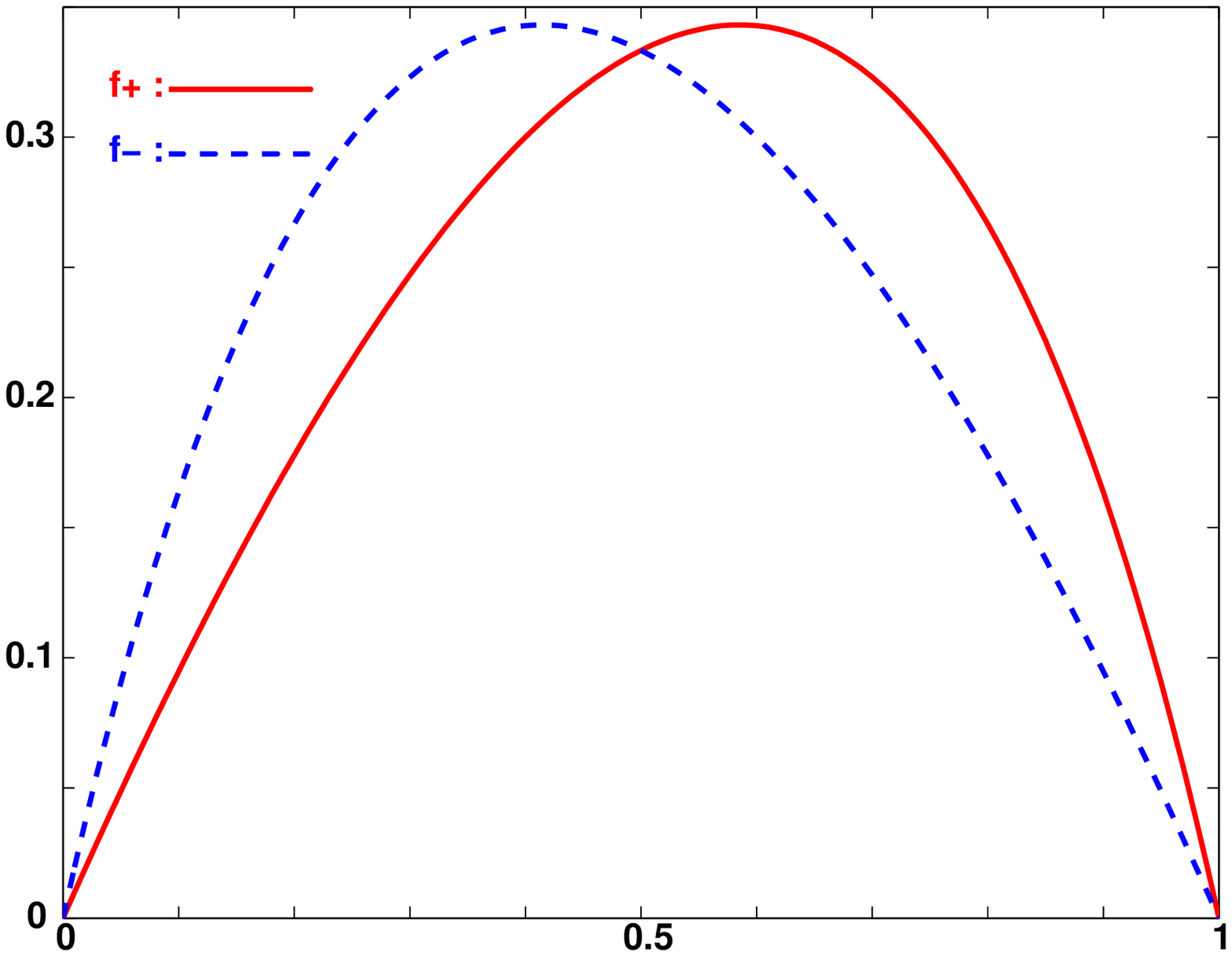}
\caption{Flux functions in experiment 2}
\end{center}
\label{e21}
\end{figure}

In this experiment, we consider the following flux functions and parameters,
\[
\begin{array}{llllllllll}
\lambda_{1}^{+}(S) &=& S && \lambda_{2}^{+}(S)&=& 2(1-S) \\
\lambda_{1}^{-}(S) &=& 2S && \lambda_{2}^{-}(S) &=& 1-S\\
g_{1} &=& 2 && g_{2}&=& 1 \\
\phi &=& 1 && q &=& 0
\end{array}
\]
In this case, we are changing the relative permeability functions across
the interface. The flux functions are shown in Fig. \ref{e21}.

\noi Observe that in this case, the flux functions intersect at the
point $0.5$ in the interior of the domain and the point of
intersection is undercompressive i.e $f^{+\prime}(0.5) > 0$ and
$f^{-\prime}(0.5) < 0$. The initial data are
$S_{0}(x) = 0.5 \quad \forall \quad x \in \R,$
so we start with a state where the light and heavy phases are fully
mixed. In this case, the entropy solution of
\cite{AJG1} is given by the constant state $0.5$ connected to the left
trace $0.42$ by a rarefaction fan on the left and the constant state
$0.5$ connected to the right trace $0.58$ on the right. Observe that
this solution satisfies the interface entropy condition (\ref{33}) as
$f^{-\prime}(0.42) = 0$ and $f^{+\prime}(0.58) = 0$.  

\begin{figure}[H]
\epsfysize=4.5cm
\epsffile{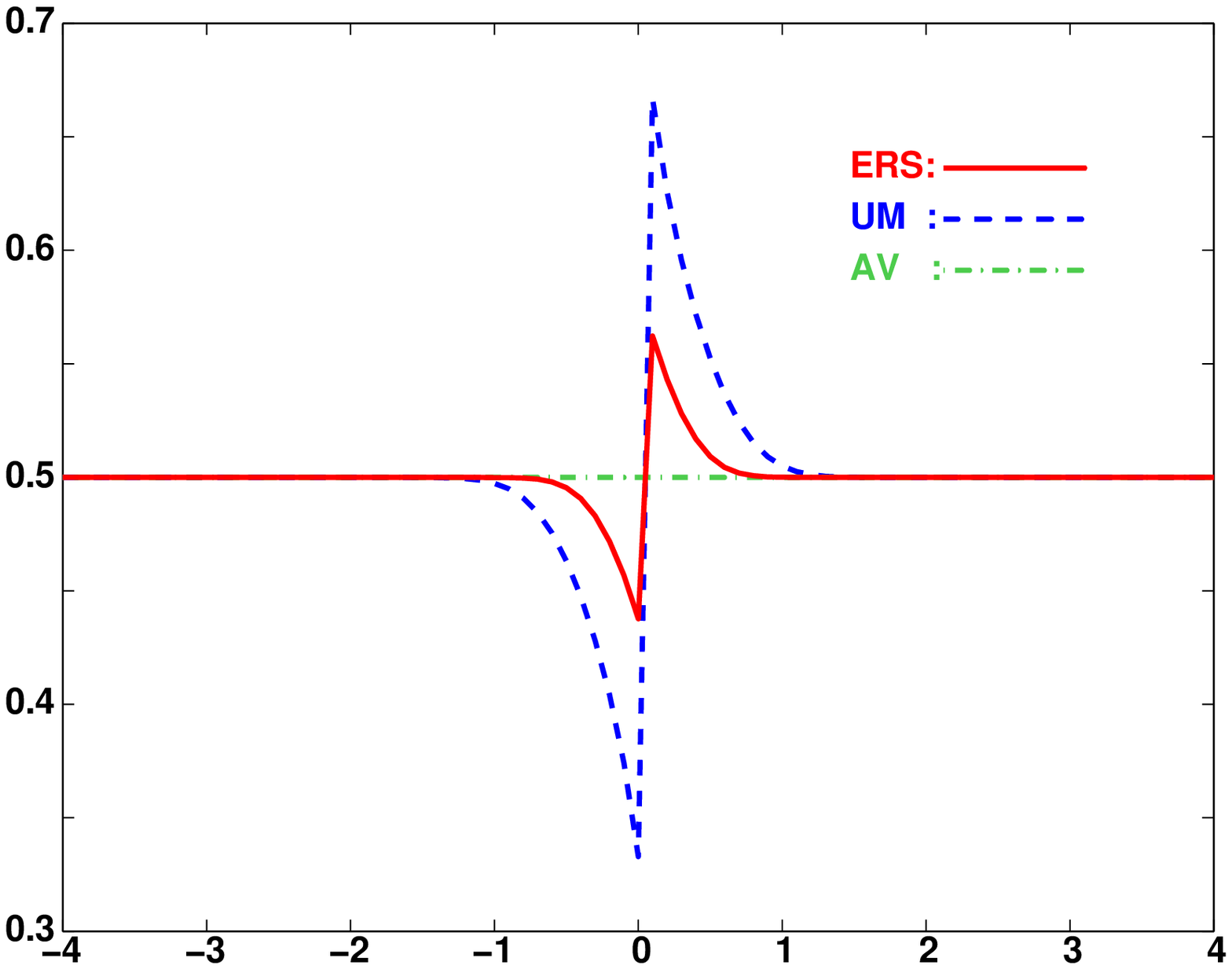} \hspace*{0.3cm} \epsfysize=4.5cm 
\epsffile{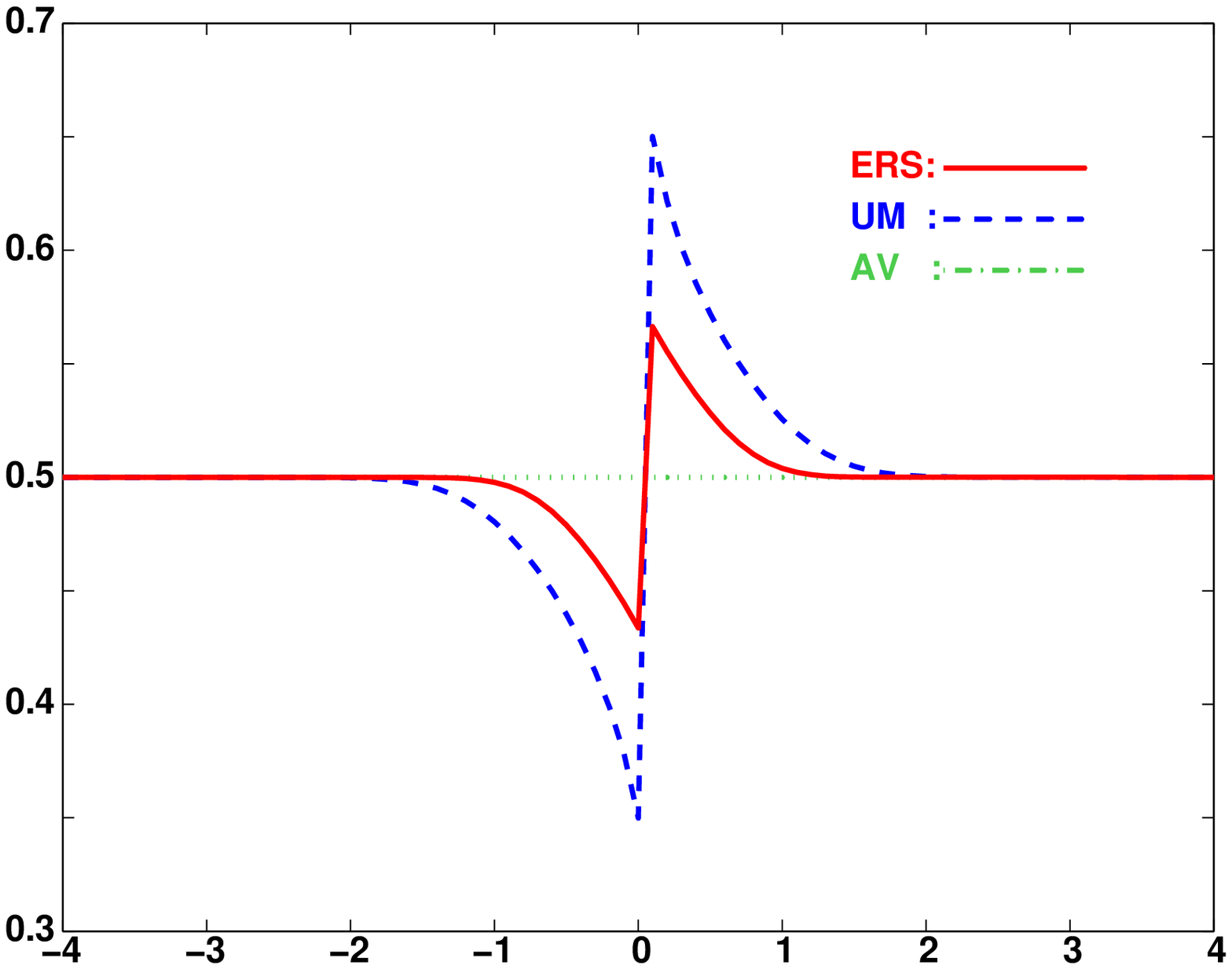}
\caption{Solutions in experiment 2 with  h=0.1 at t=1.5 and t=3}
\label{e22}
\end{figure}

\begin{figure}[H]
\epsfysize=4.5cm
\epsffile{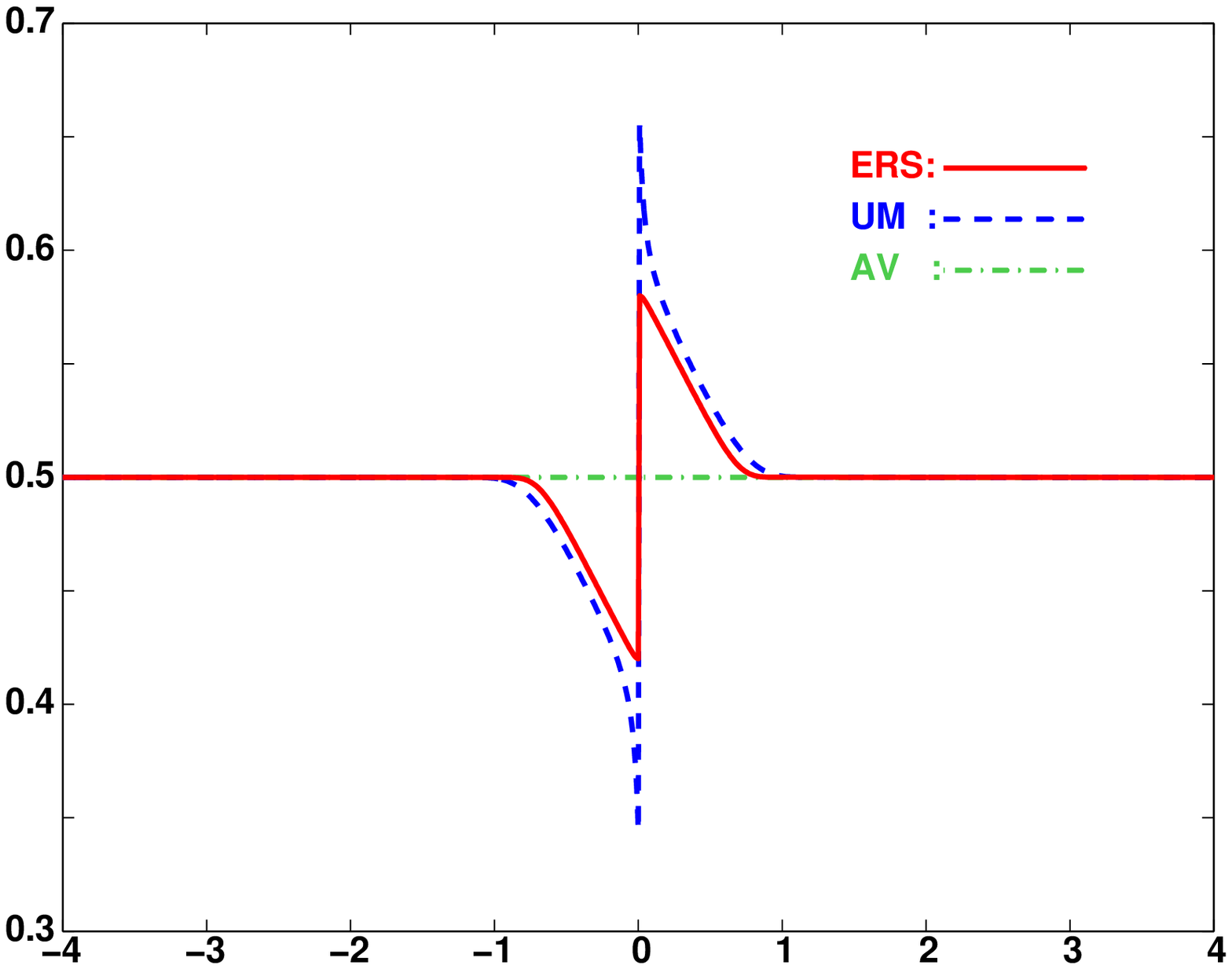} \hspace*{0.3cm} \epsfysize=4.5cm 
\epsffile{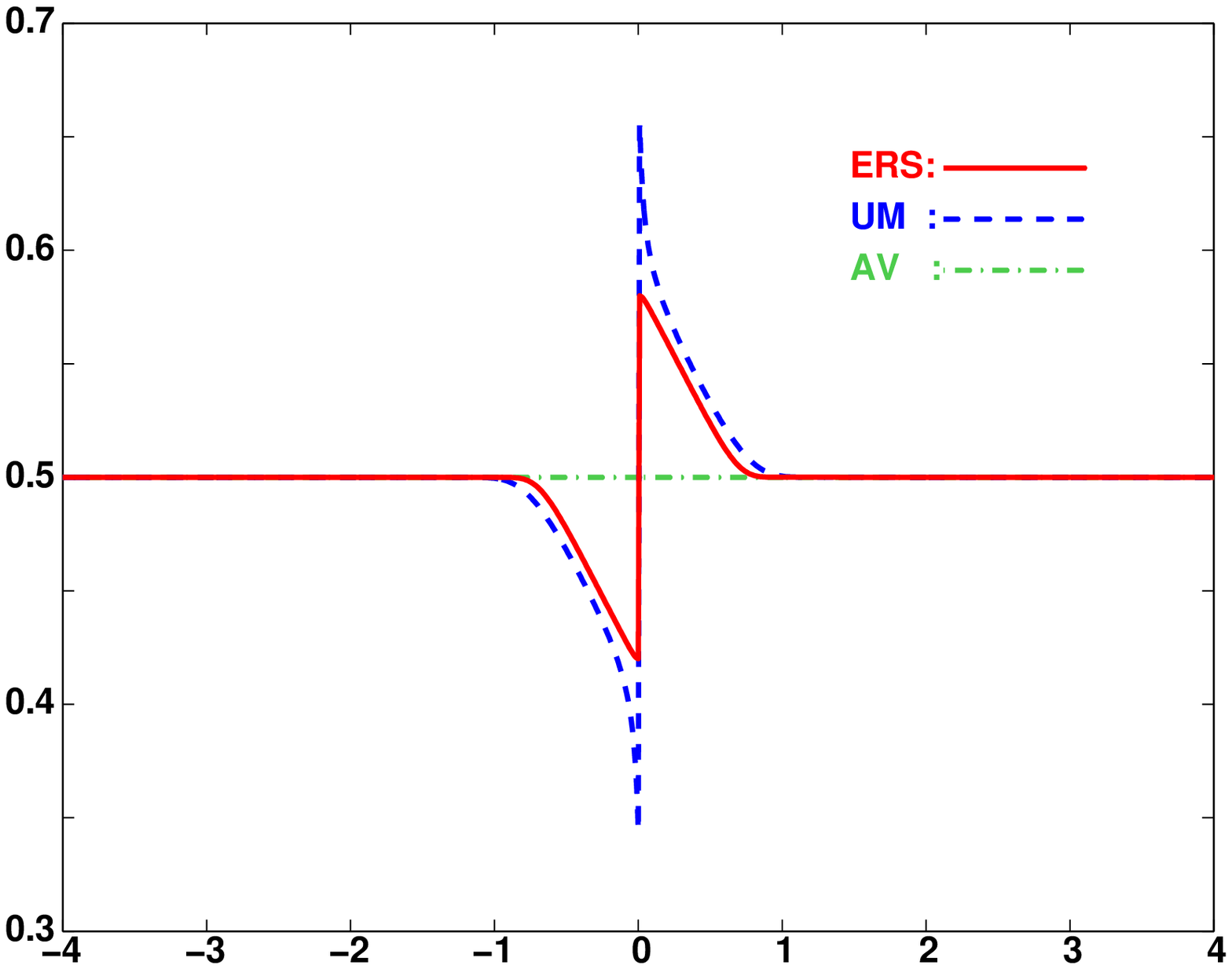}
\caption{Solutions in experiment 2 with h = 0.01 at t=1.5 and t=3}
\label{e23}
\end{figure}

As the flux functions satisfy the ``crossing'' condition of
\cite{KRT2}, we can apply the Kruzkhov type condition of \cite{KRT2}
to get that their entropy solution in this case is given by $S \equiv
0.5$. This implies that there is no flow in the medium which is
unnatural as noticed in \cite{KASS1}. This is one situation where the
above entropy theories differ and the entropy theory of \cite{AJG1}
captures the physically relevant solution. We have computed
the solutions using all the three schemes to obtain the results as
shown in Figures \ref{e22} and \ref{e23}. Fig. \ref{e22})
shows the solutions obtained by schemes ERS, UM and AV with $h=0.1$ and
the CFL parameter $\lambda=0.125$. We show the computed solutions at
times $t=1.5$ and $t=3$ respectively. As can be observed in
Fig. \ref{e22}, the solution obtained by AV is the constant state
$0.5$ in accordance with the entropy theory of \cite{KRT2}. The
solution computed by ERS converges towards the entropy solution as
discussed above with a good resolution of the interface discontinuity
and some numerical diffusion at the rarefactions. The solution
obtained by UM shows the same qualitative behaviour as that
calculated by ERS although the left
trace is $0.35$ which is well below the left trace of the solution i.e
$0.42$. Similarly the right trace of the UM solution is $0.65$ which
is above the required right trace of $0.58$. This is again indicating
the evidence for UM of a numerical boundary layer at the interface which was
noticed in experiment 1. 

In order to get a better estimate of the boundary layer, we shrink the
mesh size to $h=0.01$ and present the results in Fig. \ref{e23}.


Again, the solution obtained by AV is the constant state $0.5$. As
expected given the convergence results presented before, the solution
obtained by ERS is almost exact. Notice that the left and right traces
at $x=0$ are exactly $0.42$ and $0.58$ as in the exact solution of
this Riemann problem showing the high resolution of the scheme. The
qualitative behaviour of the solution obtained by UM is again
similar to that of ERS. But the boundary layer at the interface
remains as the left trace is still below $0.35$ (well below $0.42$)
and the right trace is still above $0.65$ well above the required
right trace of $0.58$. This suggests to us that the boundary layer
remains intact as $h \rightarrow 0$ and the traces at the interface
are different from the expected traces, although, the width of this
boundary layer shrinks with a reduction in the mesh size. This
suggests that the limit solution obtained by UM converges to the
entropy solution of \cite{AJG1} in an integral sense.\\

\begin{figure}[H]
\begin{center}
\epsfysize=6cm
\epsffile{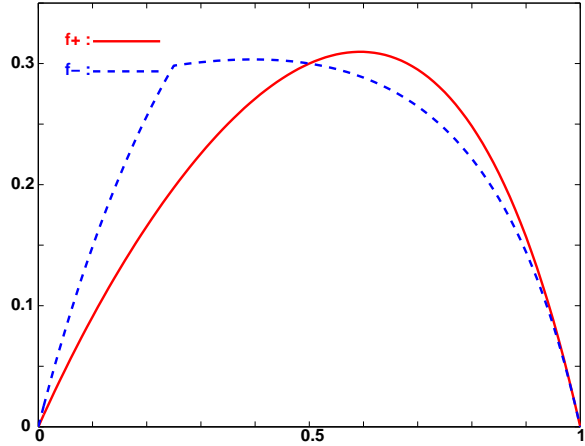}
\caption{Flux functions in experiment 3}
\end{center}
\label{e31}
\end{figure}

In experiment 2, we considered flux functions where the solutions
obtained by UM converged to the entropy solution of \cite{AJG1} in
an integral sense. The crucial point of the previous experiment was
that the fluxes intersect in the interior of the interval (0,1) and the point of
intersection was undercompressive. Next, we consider a situation of
the similar type where solutions computed by UM seem to behave very
differently.\\ 
\noi {\bf Experiment 3}
In this experiment we consider the flux functions and parameters given by,
\[
\begin{array}{llllllllll}
\lambda_{1}^{+}(S) &=& S && \lambda_{2}^{+}(S)&=&(1-S^{2}) \\
g_{1} &=& 2 && g_{2}&=& 1 \\
\phi &=& 1 && q &=& 0
\end{array}
\]
\[
\begin{array}{lllllllll}
\lambda_{1}^{-}(S)&=& 1.75S &{\rm if}& S \leq 0.25 \\
&=&0.25S + 0.375 &{\rm if}& S \geq 0.25\\
\lambda_{2}^{-}&=&1-S^{2}
\end{array}
\]
The flux functions are schown in Fig. \ref{e31}. Notice that in this
case, $f^{-}$ and $f^{+}$ intersect at $0.5$ and the intersection is
undercompressive. This is a situation which looks similar to that of the previous
numerical experiment. 
Again we start with the constant initial data given by
$S_{0}(x) = 0.5 \quad \forall \quad x \in \R$.

\begin{figure}[H]
\epsfysize=4.5cm
\epsffile{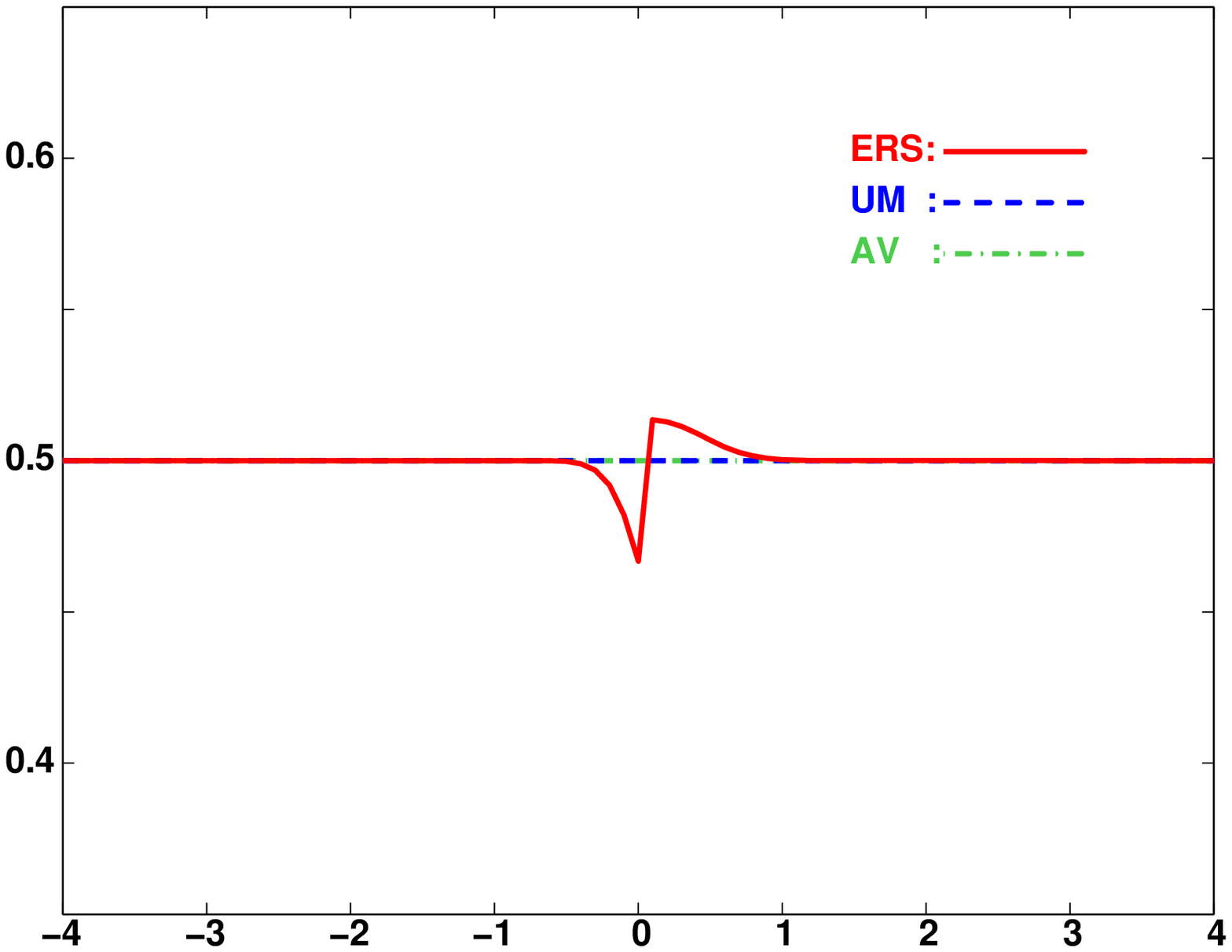} \hspace*{0.3cm} \epsfysize=4.5cm 
\epsffile{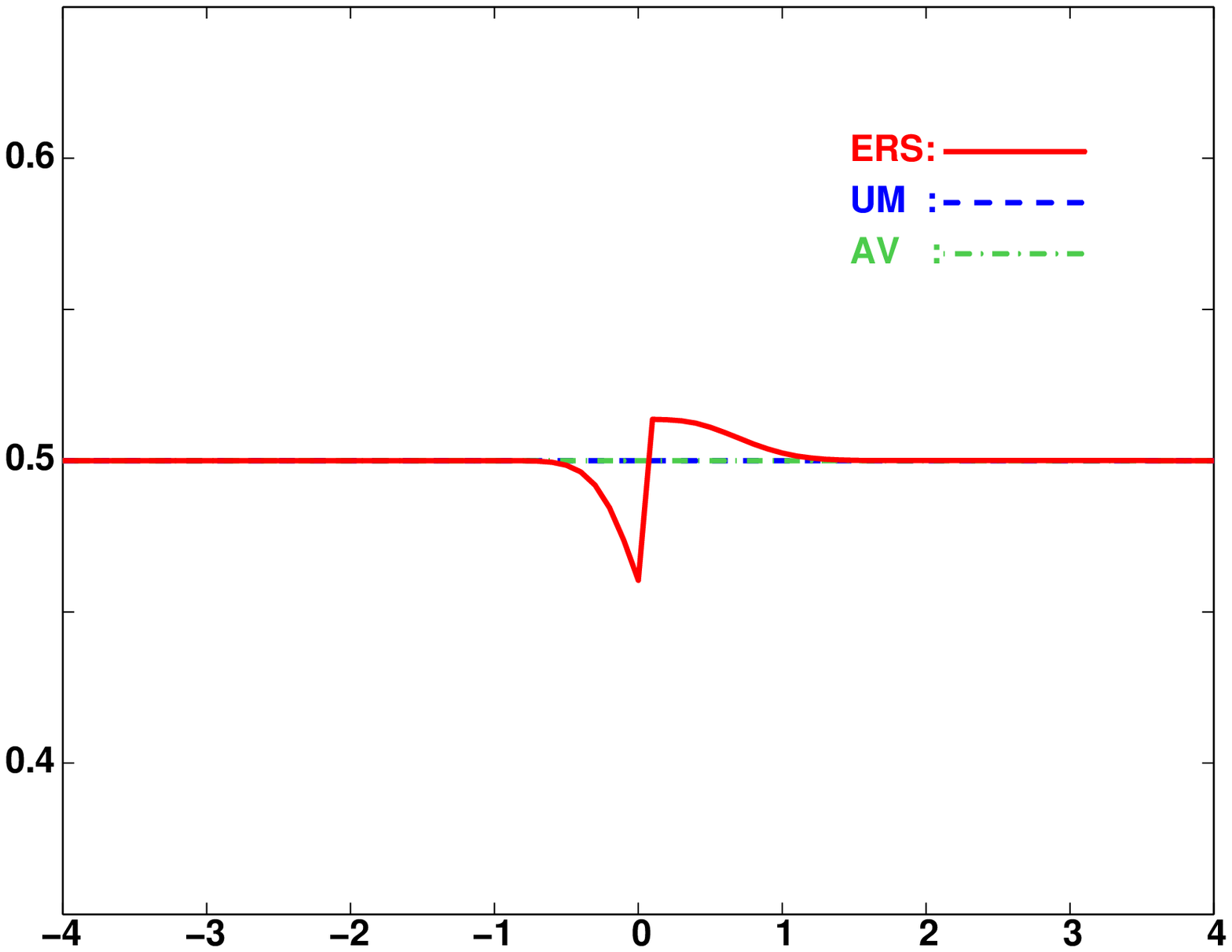}
\caption{Solutions in experiment 3 with h = 0.1 at t=2.5 and t=3.75}
\label{e32}
\end{figure}

As in experiment 2, the entropy solution of \cite{AJG1} consists of
the constant state $0.5$ connected to the left trace $0.45$ by a
rarefaction on the left, a steady discontinuity at $x=0$ connecting
the left trace $0.45$ to the right trace $0.54$,and the constant state
$0.54$ being connected by another rarefaction to the constant state
$0.5$ on the right. As the ``crossing condition'' is satisfied, the
entropy solution of \cite{KRT2} is just the constant $S \equiv
0.5$. We show  in figure (\ref{e32}) the results obtained by all three schemes with
$h=0.1$ and the CFL $\lambda = 0.125$. 

As noticed in Fig. \ref{e32}, the solution computed by ERS
approximates the entropy solution of \cite{AJG1} while the solution
computed by AV is the constant $0.5$. But what is really surprising is
that the solution computed by UM is also the constant $0.5$. In fact,
this example has been constructed in such a way that
$\lambda_{1}^{-}(0.5)=\lambda_{1}^{+}(0.5)$ and
$\lambda_{2}^{-}(0.5)=\lambda_{2}^{+}(0.5)$. Hence from the very
definition of the upstream mobility flux, it is easy to check that the
solution computed by UM remains the constant $0.5$ at all time
steps. 

Thus so far we have shown two experiments involving fluxes with an
undercompressive intersection in which the entropy solutions of
\cite{AJG1} and \cite{KRT2} differ. In experiment 2, the solutions
computed by UM flux seems to converge to the entropy solution of
\cite{AJG1} in an integral sense, though not pointwise whereas, in
experiment 3, the UM flux gives the constant solution which has been
considered in literature as unphysical (see \cite{KASS1}) and
converges to the entropy solution of \cite{KRT2}. Despite similar flux
geometry, this inconsistent behaviour of the UM flux indicates the
difficulties of characterizing the limit solutions computed by the
scheme.\\

The above numerical experiments clearly show that the inconsistent
behaviour of the UM flux when the fluxes intersect in the interior of the
interval (0,1) and the
point of intersection is undercompressive. 
We now investigate another
type of flux geometry in which the flux functions intersect and the
point of intersection is overcompressive. In this case, the limit
solution obtained with the UM flux also shows an inconsistent entropy behaviour.\\ 
\noi {\bf Experiment 4}
In this experiment, we consider the following flux functions and parameters,
\[
\begin{array}{llllllllll}
\lambda_{1}^{+}(S) &=& 2S && \lambda_{2}^{+}(S)&=& (1-S) \\
\lambda_{1}^{-}(S) &=& S && \lambda_{2}^{-}(S) &=& 2(1-S)\\
g_{1} &=& 2 && g_{2}&=& 1 \\
\phi &=& 1 && q &=& 0
\end{array}
\]
\begin{figure}[H]
\begin{center}
\epsfysize=6cm
\epsffile{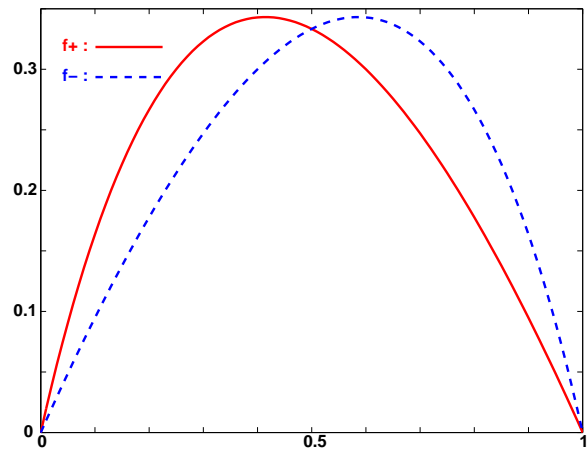}
\caption{Flux functions in experiment 4}
.eps \end{center}
\label{e41}
\end{figure}
The flux functions are shown in Fig. \ref{e41}. Observe that $f^{-}$
and $f^{+}$ intersect at $0.5$ and that the intersection is overcompressive
i.e $f^{-\prime}(0.5) > 0$ and $f^{+\prime}(0.5) < 0$. We consider the
initial data $ 
S_{0}(x) =  \left\{ \begin{array}{lcl} 2/3 &{\rm if}& x < 0\\
1/3 &{\rm if}& x > 0.
\end{array} \right. $

\begin{figure}[H]
\epsfysize=4.5cm
\epsffile{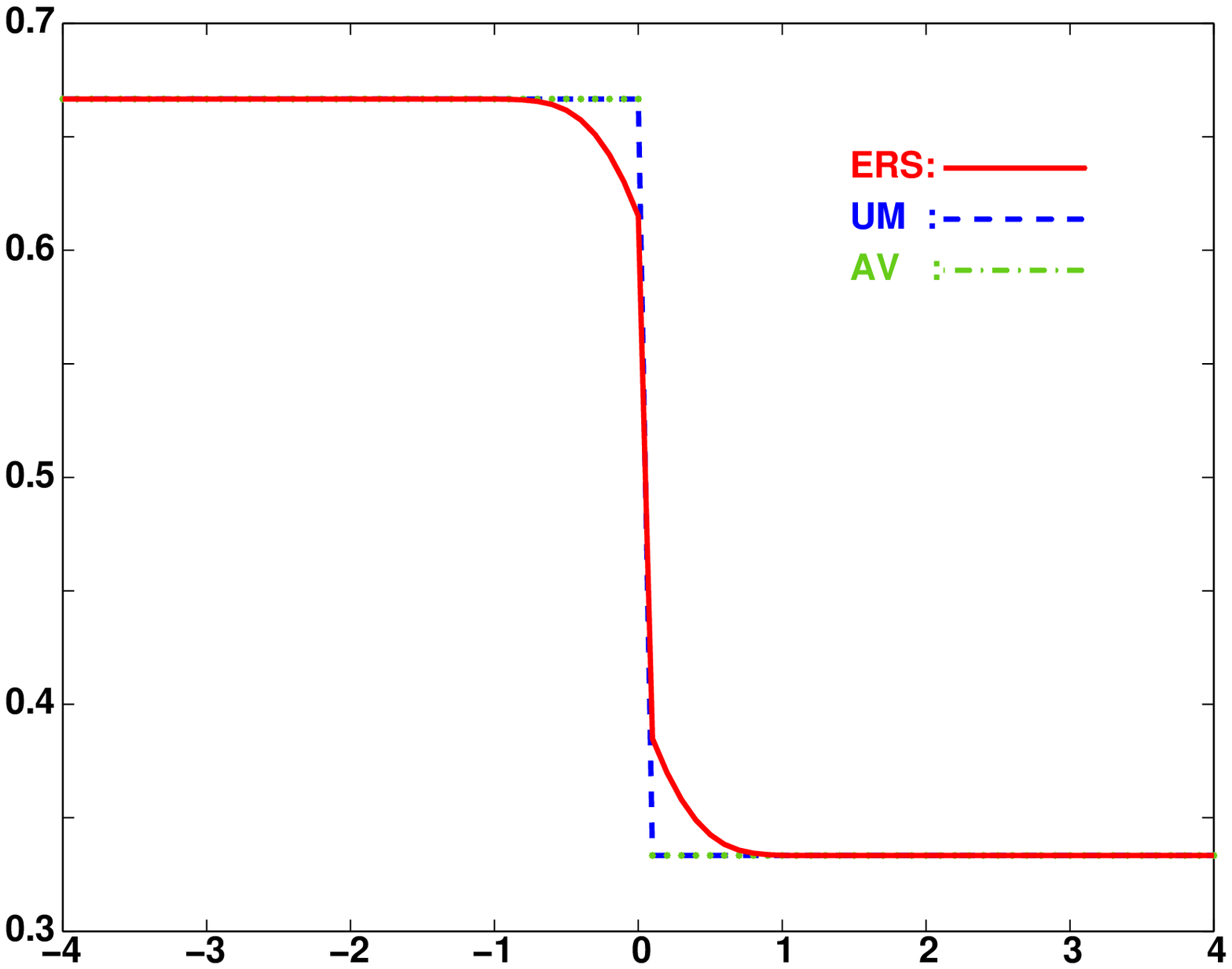} \hspace*{0.3cm} \epsfysize=4.5cm 
\epsffile{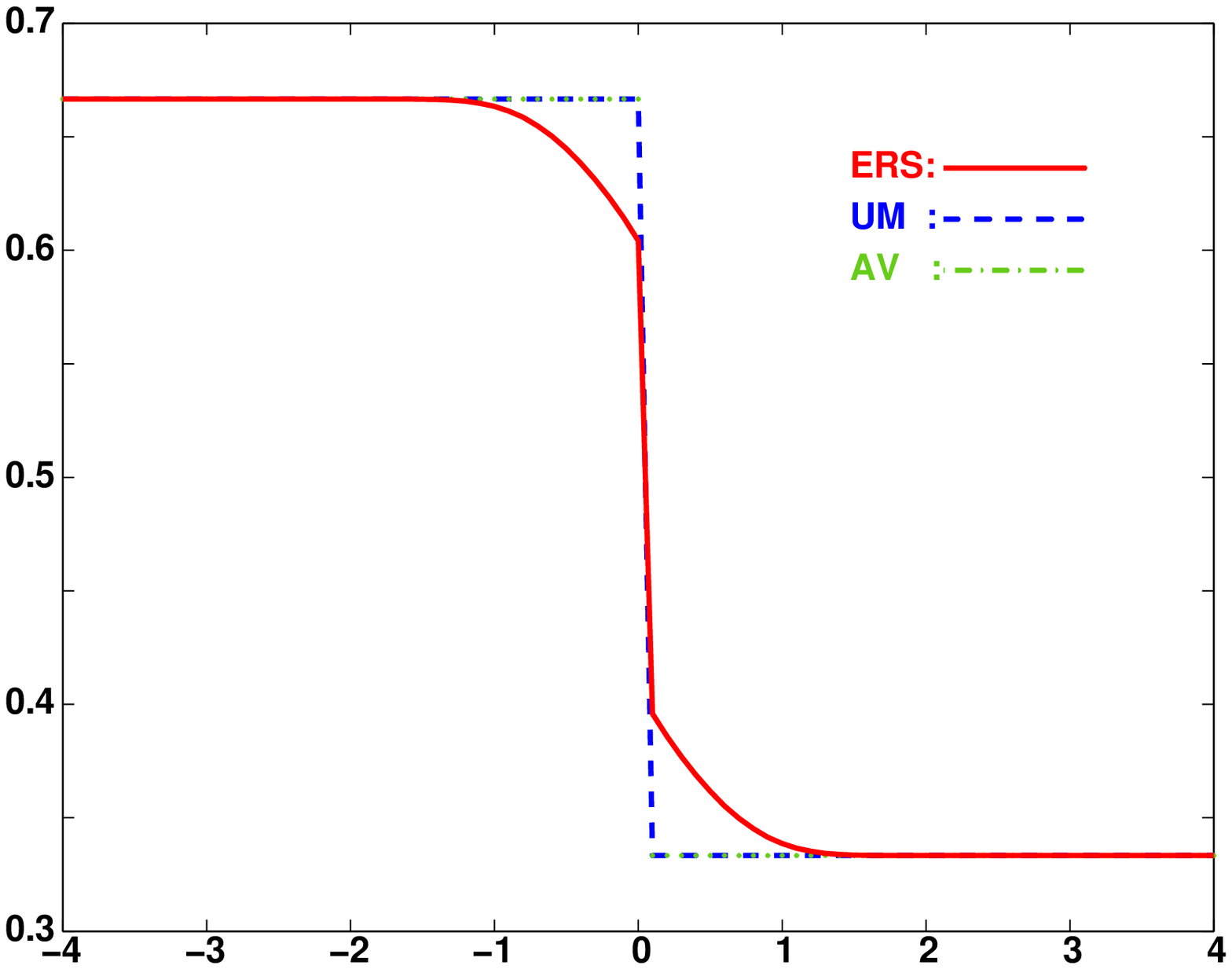}
\caption{Solutions in experiment 4 with h = 0.1 at t=1.5 and t=3}
\label{e42}
\end{figure}

The entropy solution of \cite{AJG1} in this case consists of the
constant state $0.66$ connected by a rarefaction to the left trace
$0.58$ on the left , a steady discontinuity at the interface between
the left trace $0.58$ and the right trace $0.42$ and the constant
state $0.66$ connected to the right trace $0.42$ on the right. Note
that the solution is not undercompressive as
$f^{-\prime}(0.58)=f^{+\prime}(0.42) \equiv 0$. We remark that the
above fluxes $f^{-}$ and $f^{+}$ do not satisfy the ``crossing
condition'' of \cite{KRT2} and the entropy theory developed in the
above reference does not apply to this situation. But we can still
compute the solutions given by AV as the scheme is well defined. We
present the solutions in figure (\ref{e42}). We consider the mesh size
$h=0.1$ and the CFL parameter is $\lambda = 0.125$. 

As shown in Fig. \ref{e42}), the solution obtained with ERS approximates
the entropy solution of \cite{AJG1}. Note that the left
and right traces are very close to the expected values of $0.58$ and
$0.42$. On the other hand, the solution computed by both UM and AV
is the steady state $2/3$ on the left and $1/3$ on the right which is
very different from that of the solution given by ERS. Observe that
this solution is undercompressive i.e $f^{-\prime}(2/3) < 0$ and
$f^{+\prime}(1/3) > 0$. The entropy theory of \cite{AJG1} avoids
solutions of this type. Also this solution differs from the solution
of the Riemann problem constructed by Diehl in \cite{D1} which in this
case is identical to the solution computed by ERS. We believe that
this undercompressive solution is unphysical and the right solution is
computed by ERS. \\

\begin{figure}[H]
\begin{center}
\epsfysize=6cm
\epsffile{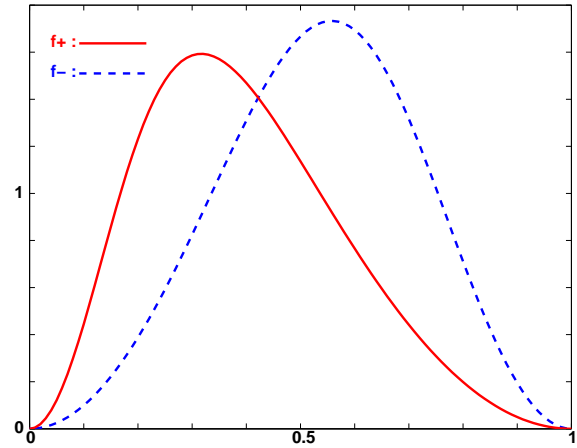}
\caption{Flux functions in experiment 5}
\end{center}
\label{e51}
\end{figure}

It is easy to show by using that $\lambda_{1}^{+}(2/3) =
\lambda_{1}^{-}(1/3)$  and $\lambda_{2}^{+}(1/3) =
\lambda_{1}^{-}(2/3)$ and the explicit definition of UM that the
solution computed by UM for all $h$ in this case is the steady state
with $2/3$ on the left and $1/3$ on the right. The natural question
that arises is whether the solutions computed by UM agree with that
of AV in the case where the flux functions intersect in an
overcompressive manner. The answer to this question is contained in
the next experiment.\\ 
\noi {\bf Experiment 5}
In this experiment, we consider the following flux functions and parameters,
\[
\begin{array}{llllllllll}
\lambda_{1}^{+}(S) &=& 50S^{2} && \lambda_{2}^{+}(S)&=& 5(1-S)^{2} \\
\lambda_{1}^{-}(S) &=& 10S^{2} && \lambda_{2}^{-}(S) &=& 20(1-S)^{2}\\
g_{1} &=& 2 && g_{2}&=& 1 \\
\phi &=& 1 && q &=& 0
\end{array}
\]
The flux functions are shown in Fig. \ref{e51}. Notice that in this
case, the flux functions intersect in the interior of the domain at
$0.46$ and the point of intersection is overcompressive. 
We consider the following initial data $
S_{0}(x) = \left\{ \begin{array}{lcl} 0.8 &{\rm if}& x < 0\\
0.2 &{\rm if}& x > 0.
\end{array} \right.$

\begin{figure}[H]
\epsfysize=4.5cm
\epsffile{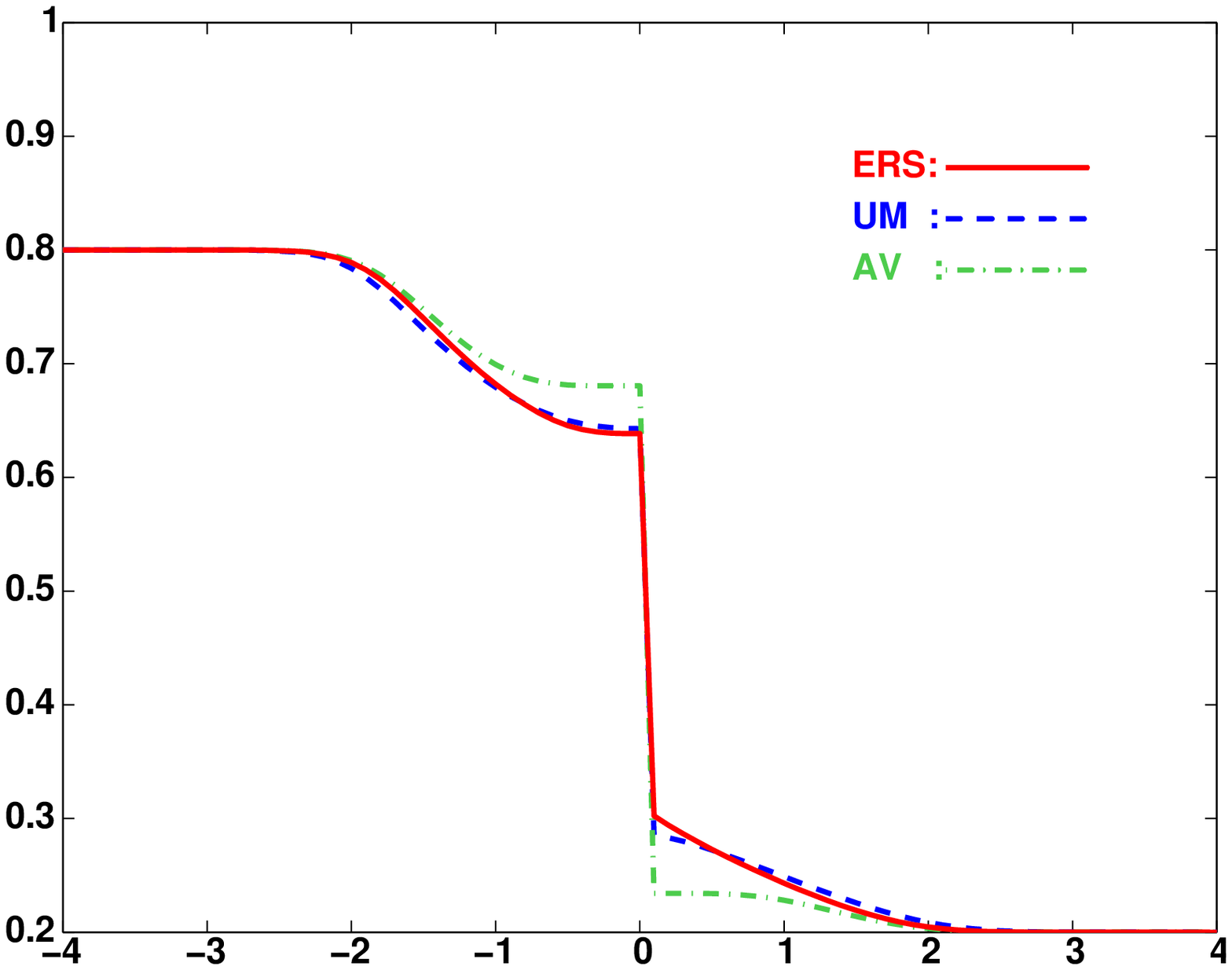} \hspace*{0.3cm} \epsfysize=4.5cm 
\epsffile{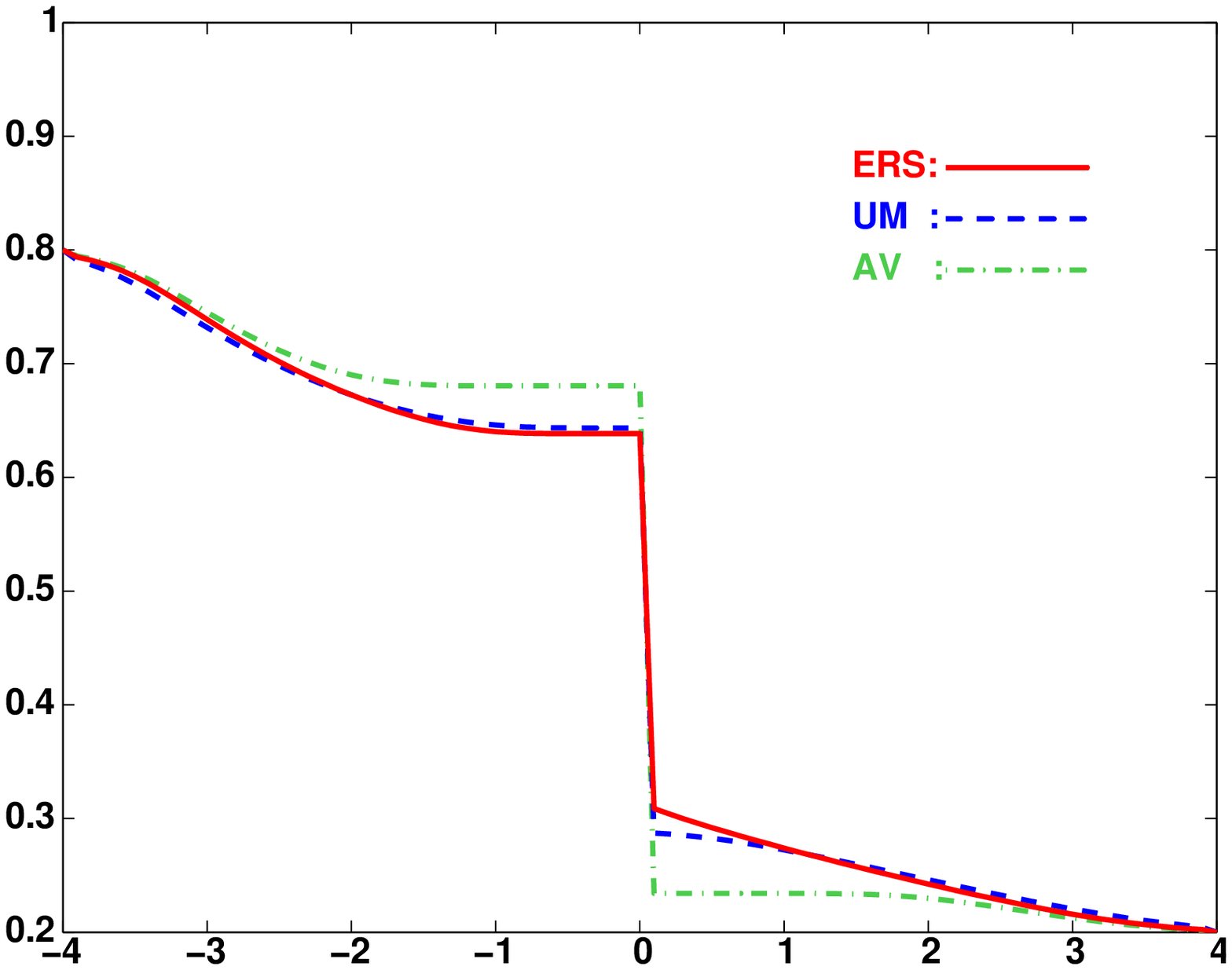}
\caption{Solutions in experiment 5 with h = 0.1 at t=0.25 and t=0.5}
\label{e52}
\end{figure}

\begin{figure}[H]
\epsfysize=4.5cm
\epsffile{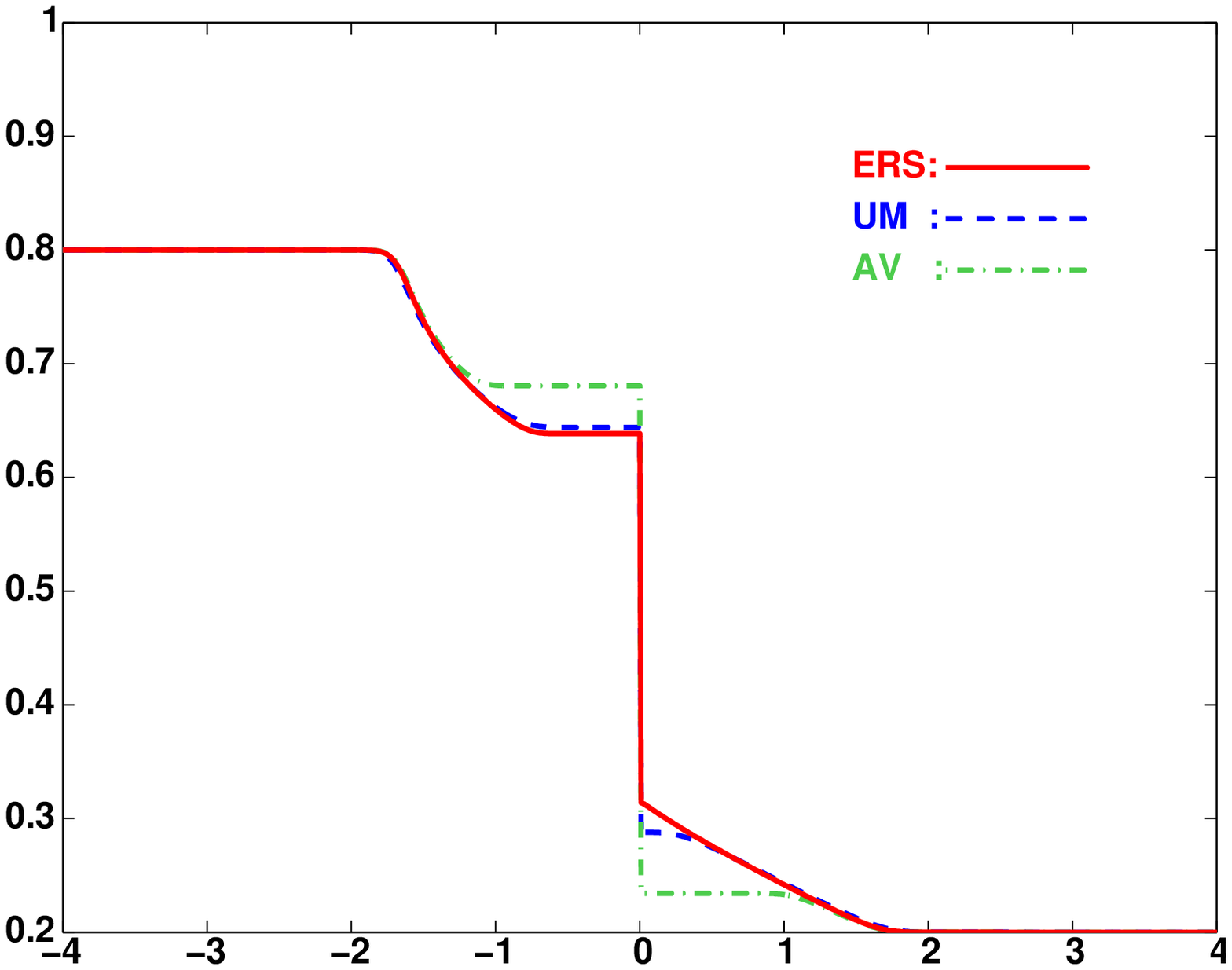} \hspace*{0.3cm} \epsfysize=4.5cm 
\epsffile{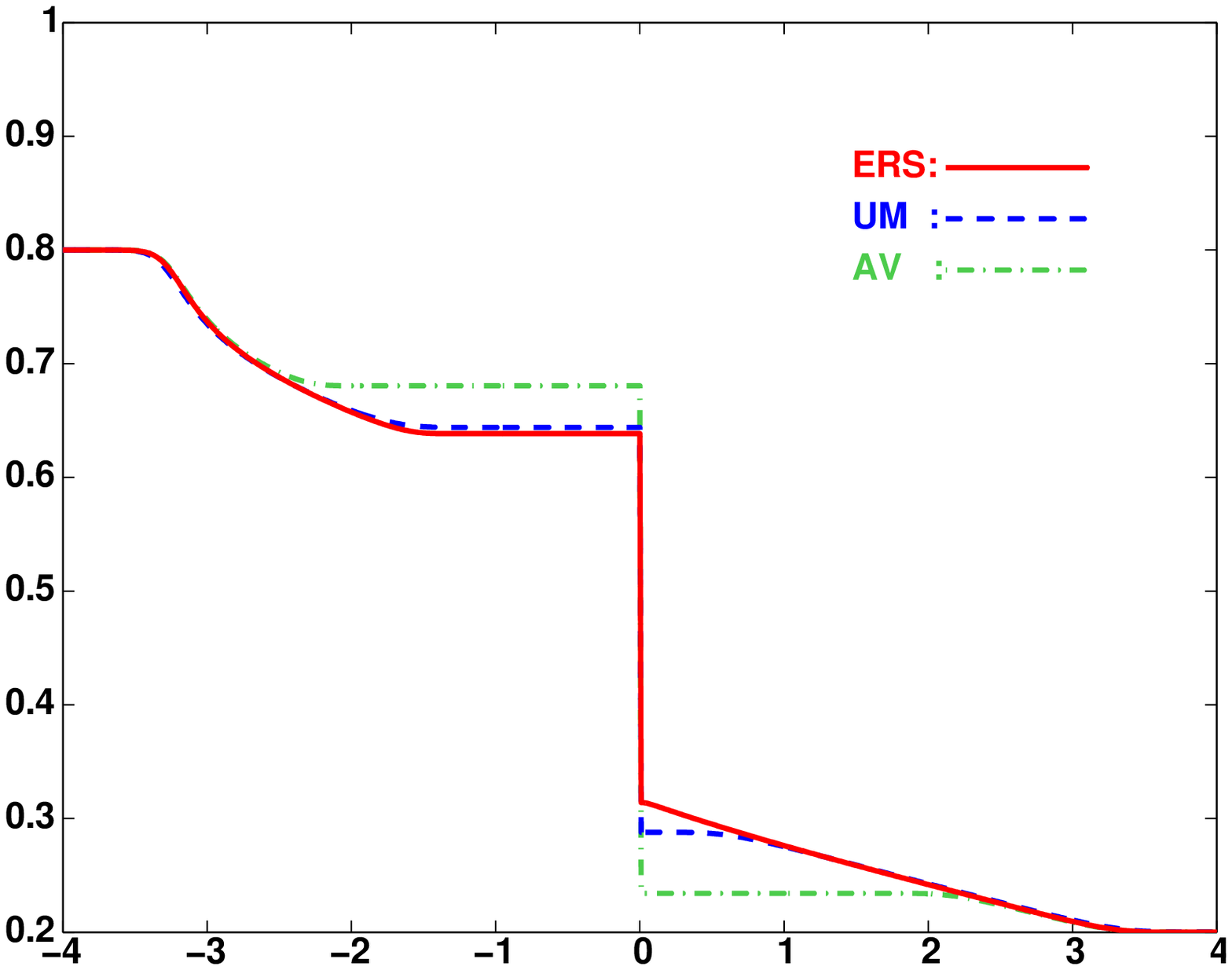}
\caption{Solutions in experiment 5 with h = 0.01 at t=0.25 and t=0.5}
\label{e53}
\end{figure}

In this case, the entropy solution of \cite{AJG1} consists of a
rarefaction joining the constant state of $0.8$ with the left trace of
$0.6$, followed by a constant state of $0.6$, a steady discontinuity
joining the left trace of $0.6$ and the right trace $0.32$ and a
rarefaction joining the right trace to that of the constant state of
$0.2$. Check that this solution is not undercompressive. The solutions
obtained by all the three schemes with $h=0.1$ and $\lambda = 1/32$
are shown in Fig. \ref{e52}. 
The solution given by the ERS flux approximates well
the entropy solution, even with a large mesh size. The solution given
by thg AV flux is quite different in this case and note that the traces
$(0.7,0.22)$ are undercompressive. On the other hand, the solutions
obtained by the UM flux are very close to those of the ERS flux besides a boundary
layer on the right. A further reduction in mesh size to $h=0.01$ shows
that the boundary layer on the right remains and the traces given by
the AV flux are undercompressive as shown in Fig. (\ref{e53}). 

To sum up about these experiments we observed the following behaviour
across the interface:
\begin{enumerate}
\item In some experiments (experiments 1,2,5) the upstream mobility
  flux may produce unphysical boundary layers and travelling waves. The traveling
  wave and the width of the boundary layer vanishes when $h
  \rightarrow 0$, while the heigth of the boundary layer may remain
  significant. Despite of these numerical artefacts the solution given
  by the upstream mobility flux remain close to that given by the
  ERS flux. This suggests that in these experiments, the
  solution calculated with the UM flux, even though it
  does not satisfy the pointwise entropy condition (\ref{33}), may
  satisfy some integral form of it. For the average flux, depending on
  the experiment, it behaves like the UM flux
  (experiments 1, 5) or it misses the interface discontinuity (experiment 2).
\item In other experiments (experiment 3, 4) the UM
  flux as well as the AV flux produces unphysical
  undercompressive solutions (experiment 4) and even misses the
  interface discontinuity (experiment 3).
\end{enumerate}
 
\section{Conclusion}
In this paper we analyzed the upstream mobility numerical flux for a finite
difference scheme when a two-phase flow crosses the interface between
two rock types. This results in a discontinuity in the flux function
with respect to the space variable. We were able to prove convergence
to a weak solution but numerical experiments show that it does not
satisfy the entropy condition of \cite{AJG1}. 

Most often the solution
given by the upstream mobility flux is close to that given by the
extended Godunov flux but numerical artefacts like boundary layers or
traveling waves perturb the solution. There are even cases when the
upstream mobility flux misses the discontinuity at the interface. The
solution given by the averaged flux is not doing any better.

\section*{Acknowledgements}
The authors would 
like to thank Professor Adimurthi and Professor G.D.Veerappa Gowda for
their useful suggestions and discussions.

\end{document}